\documentclass[a4paper,11pt]{article}

\usepackage{color}
\usepackage{graphicx}
\usepackage{amsfonts}
\usepackage{amssymb}
\usepackage{enumerate}
\usepackage{amsmath}
\usepackage{amsthm}
\usepackage[english]{babel}
\usepackage{isolatin1}

\newcommand{\RR}{\mathbb{R}}
\newcommand{\R}{\mathbb{R}}

\newcommand{\cp}{{\:\stackrel{P}{\longrightarrow}\:}}
\newcommand{\fidi}{{\:\stackrel{\mathrm{fidi}}{\longrightarrow}\:}}
\newcommand{\indi}{{1\!\!1}}
\newcommand{\rvA}{\mathbf{A}}
\newcommand{\NN}{\mathbb{N}}

\newcommand{\bz}{\mathbf{z}}
\newcommand{\bu}{\mathbf{u}}

\newcommand{\F}{{\mathcal{F}}}
\newcommand{\calF}{\mathcal{F}}

\newcommand{\calB}{\mathcal{B}}
\newcommand{\calT}{\mathcal{T}}
\newcommand{\calR}{\mathcal{R}}
\newcommand{\calC}{\mathcal{C}}
\newcommand{\calS}{\mathcal{S}}

\newcommand{\calN}{\mathcal{N}}
\newcommand{\calP}{\mathcal{P}}
\newcommand{\calE}{\mathcal{E}}
\newcommand{\G}{\mathcal{G}}
\newcommand{\K}{\mathcal{K}}
\newcommand{\calG}{\mathcal{G}}
\newcommand{\calK}{\mathcal{K}}

\newcommand{\inter}[1]{\overset{\circ}{#1}}

\newcommand{\xx}{\mathbf{x}}
\newcommand{\yy}{\mathbf{y}}
\newcommand{\xun}{\mathbf{x_1}}
\newcommand{\xn}{\mathbf{x_n}}
\newcommand{\yun}{\mathbf{y_1}}
\newcommand{\yn}{\mathbf{y_n}}

\newcommand{\Comp}[1]{#1^{\mbox{c}}}

\def\Card{\hbox{Card}}

\def\cqfd{\qquad {\ \vbox{\hrule\hbox{%
   \vrule height1.3ex\hskip0.8ex\vrule}\hrule
  }}\par}

\newcounter{hyp}
\newenvironment{hyp}[1]{\refstepcounter{hyp}\begin{itemize}\item[{\bf
      (A-\arabic{hyp})}] \label{hyp:#1}}{\end{itemize}}
\newcommand{\refhyp}[1]{{(A-\ref{hyp:#1})}}
\newtheorem{theorem}{Theorem}
\newtheorem{proposition}{Proposition}
\newtheorem{definition}{Definition}
\newtheorem{lemma}{Lemma}
\newtheorem{corollary}{Corollary}
\theoremstyle{remark}
\newtheorem{remark}{Remark}

\begin{document} 
\title{The dead leaves model : \\
general results and limits at small scales \\}
\author{Yann Gousseau, Fran\c{c}ois Roueff}
\maketitle
\begin{center}
TSI / CNRS URA 820 \\
Ecole Nationale Sup\'erieure des T\'el\'ecommunications \\
46 rue Barrault, 75634, Paris Cedex 13, France
\end{center}

\baselineskip 16pt

\setcounter{tocdepth}{2}
\tableofcontents

\newpage

\noindent {\bf Abstract} In this work, we introduce a random field in view of natural image
  modeling, obtained as a limit of sequences of dead leaves
  models, when considering arbitrarily small or big
  objects. The dead leaves model, introduced by the Mathematical Morphology school, consists in the superposition of random closed
  sets, and enables to model the occlusion phenomena. When combined with
  specific sizes distributions for objects, they are known to provide adequate models for
  natural images. However this framework yields a small scales cutoff and a limit random field is introduced by letting this cutoff tend to zero. We first give
  a rigorous definition of the dead leaves model, and compute the
  probability that $n$ compacts are included in distinct visible parts,
  which characterizes the model. Then, we derive our limit model and some of its property,
  and study its regularity.

\section{Introduction and motivations}

The structure of natural images is very specific, and strongly differs from the
one of speech signals for instance. Most statistics of natural images exhibit
non-gaussianity, as well as scaling properties. These two phenomena may
for instance be easily observed on the distribution of the gradient
of images gray levels (\cite{Ruderman94}, \cite{turiel98}). Other quantity
bearing these properties include the power spectrum (\cite{Kretzmer52}, \cite{Field87}), wavelet
coefficients (\cite{Simoncelli97}, \cite{Mumford98}), morphological
quantities (\cite{Gousseau99}) or the distribution of local patches
(\cite{GemanKoloydenko}). Non-gaussiannity is strongly related to the occlusion
phenomenon. Indeed, in the process of image formation, objects hide themselves
depending on where they lie with respect to the camera, which differs totally
from an additive generation. This phenomenon leads to peculiar twodimensional
structures such as homogeneous regions, borders and T-junctions. Besides, the scaling properties of an image may also be
seen as a result of scaling properties present in nature.

Several study (\cite{Ruderman96}, \cite{Gousseau99}, \cite{MumfordLee}) show (either theoretically or experimentally) that most of
natural images statistics may be reproduced though the use of a simple model
of images, consisting in the sequential superposition of random objects, the
dead leaves model of Mathematical Morphology. The mere nature of the model
enables the reproduction of characteristic structures of natural images
(onedimensional discontinuities, homogeneous zones). Moreover, the use of a power law
distribution $x^{-\alpha}$ for the sizes of objects enables to reproduce
scaling phenomena (note that this distribution of object sizes was also considered in
\cite{chi}, but without occlusion). In \cite{MumfordLee}, a version of the model corresponding to strict
scale invariance is considered (that is $\alpha=3$) whereas \cite{Ruderman96}
and \cite{Gousseau99}
consider $\alpha$ as a parameter of the model. In both cases, the presence of small and
large scale "cut-off" sizes is assumed, that is to say that objects sizes are
bounded. 

In this article, we present a new model for natural images, that is obtained
from a dead leaves model with scaling properties when letting the small scale
"cut-off" frequency tend to 0. By doing this, we model the small scales
properties in a non-trivial way, and we are then in a position to
study the regularity of images from a functional analysis point of
view. 

In the first part, we recall some results on random closed sets and random
tessellations, and then introduce the notion of a
"colored" (or "textured") tessellation, and give a characterization of the
convergence of a sequence of such random fields in the sense of finite
dimensional distributions. Then, we define the dead leaves model, originally introduced by
G. Matheron, \cite{Matheron68}, and study it in the framework of random
tessellations. In
particular we derive a new result enabling us to completely characterize the discontinuity set 
in the sense of random closed sets. In the last part, we introduce the specific dead
leaves model with scaling properties, and derive some of its property in the
case where objects are bounded. Then we look at the convergence of this model
when objects sizes tends to 0. We first show that such a limit model is not
well described by geometric means (i.e. by the theory of random closed sets)
because there are small objects everywhere. Then, in some cases, we show that
the finite dimensional distributions of the colored model
still exhibit interesting behavior, and thus study its convergence from this
point of view. Eventually, we study the regularity of this limit process using
Besov spaces.

\section{Random closed sets and random tessellations}
\label{sec:racs}

In this section, we recall some facts about random closed sets, and give
a definition of a random tessellation (that is a random covering of $\R^2$ by
disjoint sets) that enables to take into account covering by nonconnected
sets. We then define the notion of a "colored" tessellation, and
study the finitedimensional convergence of sequences of such processes.
  
\subsection{Random closed sets}
\label{racs}
We first recall the definition of random closed sets (RACS).
Let $\F$, $\G$ and $\K$ be respectively the sets of all closed, open and
compact sets of $\R^{n}$. Let us set for any $A \subset \R^{n}$,
\begin{align*}
\mathcal{F}^{A}=\{F\in\mathcal{F}\,:\,F\cap A=\emptyset\}, \\
\mathcal{F}_{A}=\{F\in\mathcal{F}\,:\,F\cap A\neq\emptyset\}.
\end{align*}
We define a topology $T_{\mathcal{F}}$ on $\F$, generated by the basis of open
sets
\[
(\mathcal{F}^{K},K\in\mathcal{K};\mathcal{F}_{G},G\in\mathcal{G}).
\]
In what follows, we
write $B_{\F}$ for the Borel $\sigma
$-algebras generated by $T_{\mathcal{F}}$.

\begin{definition}
\label{def:racs}
A random closed set of $\R^{n}$ is a
measurable function from a probability space $(\Omega,T,P)$ into $(\F,B_{\F})$.
\end{definition}
We have the following characterization of RACS (see \cite{Matheron75}):

\begin{proposition}
\label{prop:CaracRACS}
$B_{\F}=\sigma(\mathcal{F}_{G},G\in\mathcal{G})=\sigma(\mathcal{F}^{K},K\in\mathcal{K})$. 
In particular, a random closed set $F$ is thus characterized by its capacity function defined for every
compact set $K$ by $G_{F}(K)=Pr(F\cap K=\emptyset)$.
\end{proposition}

\paragraph{proof}
Using that for any $K\in\mathcal{K}$, there exists a sequence $(G_n)$ in $\calG$ such that $G_n\downarrow K$ and that,
for all $G\in\calG$, there exists a sequence $(K_n)$ in $\calK$ such that $K_n\uparrow G$, we obtain, respectively, 
for all $K\in\mathcal{K}$, $\mathcal{F}^{K}\in\sigma(\mathcal{F}_{G},G\in\mathcal{G})$ and, for all $G\in\mathcal{G}$,
$\mathcal{F}_{G}\in\sigma(\mathcal{F}^{K},K\in\mathcal{K})$. Hence the result.
\cqfd

\subsection{Tessellations}

A tessellation of $\RR^2$ is defined as follows.

\begin{definition}
\label{def:tessellation}
Let $T=\{F_i\}$ be a collection of closed sets. We say that $T$ is a
tessellation of $\RR^2$ if
\begin{enumerate}[(i)]
\item\label{tess1} $\bigcup_i F_i = \RR^2$.
\item\label{tess2} for all $i \neq j$, $\inter{F}_i \cap \inter{F}_j  =
\emptyset$, where $\inter{F}$ denotes the interior of $F$.
\item\label{tess3} $\Card\{ i \,:\, F_i \cap A \neq \emptyset \} < + \infty$
for all bounded $A \subset \RR^2$.
\end{enumerate}
Moreover we define the boundary of $T$ by 
$\partial T = \bigcup_i \partial F_i$, where $\partial F_i$ denotes the topological
boundary of $F_i$. 
\end{definition}
Observe that, by~(\ref{tess3}), the number of non-empty sets in a tessellation is necessarily countable. Hence, without loss of generality, we may index a tessellation by integers.
Let $\calT$ be the set of all tessellations indexed by $\NN$. 

\begin{lemma}
We have $\calT\in\mathcal B_{\calF^\NN}$, where $\calB_{\calF^\NN}$ is the
$\sigma$-field generated by the cylinders of the product space $\calF^{\NN}$.  
\end{lemma}
\paragraph{proof}
We notice that~((\ref{tess1}) and~(\ref{tess3})) is equivalent to  
\begin{equation}
\label{eq:equiv_tess1}
\forall N\in\NN,\exists n\in\NN, \bigcup_{i=1}^n F_i \supset B(0,N), 
\end{equation}
and
$$
\forall N\in\NN,\exists n\in\NN,\forall p\in\NN,  
\Card\{ i \in\{1,\ldots,p\}\,:\, F_i \cap B(0,N) \neq \emptyset \} 
\leq n,
$$
where $B(\xx,r)$ is the open ball centered at $\xx$ with radius $r$.  
The second condition clearly define sets $B_{\calF^\NN}$. 
Now we know from \cite[section 1-2]{Matheron75} that $(F,F')\mapsto F\cup F'$ is continuous and thus is measurable and
that $F\mapsto \Comp{(\inter{F})}$ is lower semi-continuous which implies
measurability. Now note that condition (\ref{tess2}) is equivalent to
$$F_i \cap F_j \subset \partial F_i \cup \partial F_j.$$
We have (\cite{Matheron75}) that $(F,F')\rightarrow F\cap F'$ is upper
semicontinuous and that $F\rightarrow \partial F$ is lower semicontinuous, and
thus that these two applications are measurable. Eventually, it remains to
prove that $(F,F')\rightarrow \indi_{F\subset F'}$ is measurable. Note that
$F\subset F'$ is equivalent to $F\cap F' = F$, and that the diagonal set
$\{(F,F), F\in \calF\}$ is a closed set since $\F$ is Hausdorff. This
proves the requested measurability. 
\cqfd
This enables us to define a random tessellation as follows: 
\begin{definition}
  A random tessellation of the plane is a measurable map $T$ from a probabilistic space $(\Omega, \mathcal S,
P)$ to $(\calF^\NN, \mathcal B_{\calF^\NN})$, such that almost surely $T \in \calT$. 
\end{definition}

Property (iii) implies that for any random tessellation $T$, $\partial T$ is a random closed set. 
In a random framework, there has been various definitions of random
tessellations (see e.g. \cite{Ambartzumian74}, \cite{Cowan80},
\cite{Moller89}, \cite{Stoyan95}), mostly when the $F_i$ are convex
polygons. A classical approach is to define $\partial T$ directly as a random closed set
without considering the underlying tessellation. However, it is not
always possible to recover the $F_i$'s from $\partial T$ (they may not be connected). 
Thus, we choose Definition~\ref{def:tessellation}, that provides a more informative
definition of a random tessellation than its border. This is needed in Section \ref{colored}
where we study the random function that arises when independently coloring each
$F_i$, in a case where they are generally not connected. Classical examples of random
tessellations (see the references in \cite[Chapter 10]{Stoyan95} and \cite{okabe:00}) include Poisson line
processes, Delaunay and Voronoi tessellations,  and the dead leaves model that we consider
in Section \ref{sec:deadleaves}.

\subsection{A colored tessellation process and its finite-dimensional convergence}
\label{colored}

In this section, 
we study the random process that arises when independently 
``coloring'' (or texturing) each part of a random tessellation. In particular,
we give simple criteria for the finite dimensional convergence of such a
process. First observe that any tessellation $T=\{F_i\}$ partitions $\R^2$ into $\partial
T\cup\bigcup_i(\inter{F}_i\backslash\partial{T})$. 

\begin{remark}
\label{rem:tessbis}
Note that we could
have given a more restrictive definition of a random tessellation where property
(\ref{tess3}) is replaced by 
\begin{equation}
  \label{eq:tess4}
\mbox{for all } i \neq j, {F}_i \cap \inter{F}_j  =
\emptyset.
\end{equation}
In this case, the plane is partitioned into $\partial
T\cup\bigcup_i\inter{F}_i$. Also note that it is easily seen that if we assume that the sets $F_i$
verify, 
\begin{equation}
  \label{eq:baldsets}
  \mbox{for all } i, F_i = \overline{\inter{F_i}},
\end{equation}
then (\ref{eq:tess4}) is implied by property (\ref{tess2}) (indeed if $x\in \inter{F_j}$, there exists a ball $B(x,r)\subset
\inter{F_j}$, so that if $\inter{F_i}\cap\inter{F_j}=\emptyset$, then $x
\notin \overline{\inter{F_i}}$). This gives an easy way to ensure that the
boundary of the tessellation does not interfere with the interiors of the
$F_i$, in which case some of the following definitions and results may be
slightly simplified.
\end{remark}

\begin{definition}\label{def:CRT}
  Let $T=\{F_i\}$ be a random tessellation and
  $(\{C(\xx)\,:\,\xx\in\RR^2\})$ be a real valued random field. 
For all $\xx \in \R^2\backslash {\partial{T}}$, denote by $i(\xx)$ the unique index such that
$\xx\in\inter{F}_{i(\xx)}$. 
Let
  $I=\{I(\xx)\,:\,\xx\in\R^2\}$ be the random field on $\R^2$ defined by
$$
\left\{
\begin{array}{ll}
I(\xx)=C_{i(\xx)}(\xx)& \mbox{ for all $\xx \in \R^2\backslash {\partial{T}}$},\\
I(\xx)=0 & \mbox{ for all $\xx\in{\partial}T$},
\end{array} 
\right.
$$
where $(\{C_i(\xx)\,:\,\xx\in\RR^2\})_i$ are independent copies of $C$, independent of $T$.
We call $I$ the colored tessellation process associated to $T$ and $C$.
\end{definition}

The finite-dimensional distributions of $I$ are easily obtained from the definition by conditioning on how the tessellation 
partitions the plane. The following notations allows for a formalization of this simple idea.
For a given tessellation  $T=\{F_i\}$, we denote by $\calR^{(T)}$ the equivalence relationship defined by the partition
$\left\{\partial{T},\,\inter{F_i}\backslash\partial{T}\right\}$, that is, for all $\xx,\yy\in\R^2$
$\xx\calR^{(T)}\yy$ if and only if either $\xx,\yy\in\partial{T}$ or
there exists $F_i$ such that $\xx$ and $\yy$ both are in $\inter{F_i}\backslash\partial{T}$.
For all $n\geq 1$, we denote by $\calC_n$ the set of all subsets of 
$\{1,\ldots,n\}$ including the empty set. For any set $A$, we let $\calS(A)$ denote the set of all partitions of $A$,
with the convention $\calS(\emptyset)=\emptyset$. For all $\kappa\in\calS(A)$ 
and all $a\in A$, we denote by $\kappa(a)$ the equivalence class of $a$, that is $\kappa(a)=\kappa(b)$ if and only if 
$a$ and $b$ belong to the same subset in the partition $\kappa$.
Let $\xun,\dots,\xn \in \RR^2$.
For all $A\in\calC_n$ and for all $\kappa\in\calS(A)$, 
we define the $n$-dimensional random variable 
\begin{equation}
  \label{eq:IAkap}
  I(A,\kappa,\xun,\dots,\xn):=(\tilde{C}(A,\kappa,1),\ldots,\tilde{C}(A,\kappa,n)),
\end{equation}
where  $\tilde{C}(A,\kappa,i)=0$ if $i\notin A$ and 
$\tilde{C}(A,\kappa,i)=C_{\kappa(i)}(\mathbf{x_i})$ otherwise, $\{C_j,\,j\in\kappa(A)\}$ being independent copies of
$C$. Define the random variable $(\rvA,K)$ to be the element of the finite states space   
$\{(A,\kappa),\,A\in\calC_n,\,\kappa\in\calS(A)\}$ defined by
$\rvA:=\{i\in\{1,\ldots,n\}\,:\,\mathbf{x_i}\notin\partial T\}$ and $K$ is the partition of $\rvA$ defined by
the equivalence relationship 
$\mathbf{x_i}\calR^{(T)}\mathbf{x_{j}}$ for all $i,j\in\rvA$. Further define 
\begin{multline}\label{eq:w}
w(A,\kappa,(\mathbf{x_i})_{i=1}^n):=Pr((\rvA,K)=(A,\kappa))= \\
Pr\left( 
\bigcap_{i=1}^n\left\{i\in A\Leftrightarrow\mathbf{x_i}\notin{\partial}T\right\}\cap
\bigcup_{i,j=1}^n\left\{\mathbf{x_{i}}\calR^{(T)}\mathbf{x_{j}}\Leftrightarrow \kappa(i)=\kappa(j)\right\}
\right).
\end{multline}

\begin{lemma}\label{lem:finitedistr}
For all $\xun,\dots,\xn \in \RR^2$,  
the distribution of $(I(\xun),\ldots,I(\xn))$ is the finite mixture of
$\{I(A,\kappa,\xun,\ldots,\xn),\,A\in\calC_n,\,\kappa\in\calS(A)\}$ defined with respective weights
$\{w(A,\kappa,\xun,\dots,\xn)\\,\,A\in\calC_n,\,\kappa\in\calS(A)\}$. 
\end{lemma}
\paragraph{proof}
Let $\xun,\dots,\xn \in \RR^2$.
This lemma is easily obtained by observing that conditioning $(I(\xun),\ldots,I(\xn))$ on $(\rvA,K)=(A,\kappa)$
precisely give $I(A,\kappa,\xun,\dots,\xn)$. Hence the result.
\cqfd

In cases where any point has probability zero to belong to the set $\partial{T}$, the weights~(\ref{eq:w}) are zero for
$A\neq\{1,\ldots,n\}$, which reduces the random vector $(I(\xun),\ldots,I(\xn))$ to a mixture of random variables
indexed by $\calS(\{1,\ldots,n\})$. However, even in this simpler case, 
we will not pursue in studying the geometrical structure of the obtained distribution. 
A convergence result will simply be obtained by observing that the couple $(A,K)$ defined in the proof above is a
deterministic function (which we do not precise) of a random field indexed by $(\R^2)^2$ which takes its value in
$\{0,1\}$ (see Lemma~\ref{lem:secondorder} below). Let us first introduce the following definition.   
 
\begin{definition}
Let $T=\{F_i\}$ be a random tessellation.
For all $n\geq1$ and for all $\xun,\dots,\xn \in \RR^2$, let  $R^{(n)}(\xun,\dots,\xn)$ denote the random variable which
takes value one if there exists $i$ such that the $n$ points $\xun,\dots,\xn$ are in $\inter{F_i}\backslash\partial T$
and takes value zero otherwise. We will say that $\{R^{(n)}(\xun,\dots,\xn)\,:\,(\xun,\dots,\xn)\in(\R^2)^n\}$ is the
$n^{\mbox{th}}$ order partition process. 
\end{definition}

We obtain the following result.

\begin{lemma}\label{lem:secondorder}
Let $\xun,\ldots,\xn\in\R^2$ and define the random variable $(\rvA,K)$ as in Lemma~\ref{lem:finitedistr}.
Then 
$$
(\rvA,K)\in \sigma(R^{(2)}(\xx,\yy)\,:\,\xx,\yy\in\{\xun,\ldots,\xn\}).
$$
\end{lemma}
\paragraph{proof}
Clearly, $(\rvA,K)$ is a deterministic function of the partition processes taken at finite samples in
$\xun,\ldots,\xn$, that is 
$$
(\rvA,K)\in \sigma(\left(R^{(m)}(\bz),\,\bz\in\{\xun,\ldots,\xn\}^m,\,1\leq m\leq n\}\right). 
$$
Now observe that, for all $m\geq2$ and for all $\yun,\dots,\yn \in \RR^2$, 
$$
R^{(m)}(\yun,\dots,\yn)=\prod_{i=1}^{n-1}R^{(2)}(\mathbf{y_i},\mathbf{y_{i+1}}).
$$ 
Furthermore, we have, for all $\yy\in\R^2$, $R^{(1)}(\yy)=R^{(2)}(\yy,\yy)$. The result follows. 
\cqfd

Thus it follows that the probabilities defined in~(\ref{eq:w}) may be computed
from the probability distribution of a finite sample of the process $R^{(2)}$. 
This allows for simple conditions to let a sequence $(I_j)_{j\in\NN}$ of colored tessellations converge to a limit field 
in the sense of finite-dimensional distributions. Let us recall that
$I_j\fidi I_\infty$ if, for all $n\geq1$ and for all $\xun,\ldots\xn\in\R^2$,
$(I_j(\xun),\ldots,I_j(\xn))$ converges to $(I_\infty(\xun),\ldots,I_\infty(\xn))$ in distribution. 

\begin{theorem}\label{thm:limitProc}
Consider a collection of tessellations $\{T_j\,:\,j\in\NN\}$ and, for
 all $\xun,\dots,\xn \in \RR^2$, denote by
 $R_j^{(2)}$ the second order process of $T_j$.
For all $j\in\NN$, let $(\{C_{j}(\xx)\,:\,\xx\in\RR^2\})$ be a sequence of real valued random fields, 
 independent of the tessellations $(T_j)$. Let us denote by $I_j$ the colored tessellation
 process associated to $T_j$ and $C_{j}$. Assume that  
\begin{enumerate}[{\rm (i)}]
\item \label{item:A1} there exists a random process $R_\infty^{(2)}$ in $\{0,1\}^{(\R^2)^2}$ such that
$R_j^{(2)}\fidi R_\infty^{(2)}$.  
\item \label{item:A2} there exists a random field $C_\infty$ in $\RR^{\RR^2}$ such that 
$C_{j}\fidi C_\infty$.  
\end{enumerate}
Then there exists a random field $\{I_\infty(\xx)\,:\,\xx\in\RR^2\}$ such that $I_j\fidi I_\infty$.
Furthermore the finite-dimensional distributions of $I_\infty$ only depends on those of $R_\infty^{(2)}$ and
$C_\infty$. 
\end{theorem}

\paragraph{proof}
We use the expression of the finite distributions given in Lemma~\ref{lem:finitedistr}. Then 
assumption~(\ref{item:A1}) and Lemma~\ref{lem:secondorder} imply that, for all $A\in\calC_n$ and for all 
$\kappa\in\calS(A)$, $w_j(A,\kappa,\xun,\dots,\xn)$ (defined as in~(\ref{eq:w})) has a limit
$w_\infty(A,\kappa,\xun,\dots,\xn)$ possibly computable from the finite distributions of
$R_\infty^{(2)}$. Assumption~(\ref{item:A2}) implies that, for all $A\in\calC_n$ and for all 
$\kappa\in\calS(A)$, $I_j(A,\kappa,\xun,\dots,\xn)$ (defined as in~(\ref{eq:IAkap}))
converges to the random variable  $I_\infty(A,\kappa,\xun,\dots,\xn)$ defined as in~(\ref{eq:IAkap}) using the random
field $C_\infty$. 
Then $(I_j(\xun),\ldots,I_j(\xn))$ converges in distribution to the limit mixture
$\{I_\infty(A,\kappa,\xun,\dots,\xn),\,A\in\calC_n,\,\kappa\in\calS(A)\}$ with 
weights $\{w_\infty(A,\kappa,\xun,\dots,\xn),\,A\in\calC_n,\,\kappa\in\calS_n\}$. Doing this for all $n$ and for all
$\xun,\dots,\xn \in \RR^2$, we see that all finite distributions converge and that the limit distribution only depend on 
those of $R_\infty^{(2)}$ and $C_\infty$. This gives the result. \cqfd

The two-dimensional distribution of the limit is detailed hereafter (as an application of the proof above).

\begin{proposition}\label{prop:BidimLimit}
Under the assumptions of Theorem~\ref{thm:limitProc}, the 
two-dimensional distributions of the limit process $I_\infty$ are given as follows. Let ${C'}_\infty$ be an independent
copy of $C_\infty$. For all $\xx,\yy\in\R^2$,  $(I_\infty(\xx),I_\infty(\yy))$ is a mixture of the random variables
$(C_\infty(\xx),C_\infty(\yy))$, $(C_\infty(\xx),{C'}_\infty(\yy))$, $(C_\infty(\xx),0)$, $(0,C_\infty(\yy))$, $(0,0)$
with respective weights  $Pr(R_\infty^{(2)}(\xx,\yy)=1)$,
$Pr(R_\infty^{(2)}(\xx,\yy)=0,R_\infty^{(2)}(\xx,\xx)=1,R_\infty^{(2)}(\yy,\yy)=1)$,
$Pr(R_\infty^{(2)}(\xx,\xx)=1,R_\infty^{(2)}(\yy,\yy)=0)$, $Pr(R_\infty^{(2)}(\xx,\xx)=0,R_\infty^{(2)}(\yy,\yy)=1)$ and
$Pr(R_\infty^{(2)}(\xx,\xx)=0,R_\infty^{(2)}(\yy,\yy)=0)$. 
\end{proposition}

\begin{remark}\label{rem:bidim}
If, for all $j\in\NN$,
any point of $\R^2$ has probability zero to belong to the set $\partial{T_j}$, then, we have $Pr(R_\infty^{(2)}(\xx,\xx)=0)=0$
for all $\xx\in\R^2$. Thus, for all $\xx,\yy\in\R^2$,  $(I_\infty(\xx),I_\infty(\yy))$ is a mixture of the two random variables
$(C_\infty(\xx),C_\infty(\yy))$ and $(C_\infty(\xx),{C'}_\infty(\yy))$
with respective weights  $Pr(R_\infty^{(2)}(\xx,\yy)=1)$,
$Pr(R_\infty^{(2)}(\xx,\yy)=0)$. 
\end{remark}

\section{The dead leaves model}
\label{sec:deadleaves}

In this section, we introduce the dead leaves model. This model was introduced
by G. Matheron in \cite{Matheron68}, and later presented by J. Serra in
\cite{Serra82} (see also \cite{Cowan94}, \cite{Jeulin}, \cite{Kendall99}). 

\subsection{Definition}

Before proceeding with the definition of the model, we recall some notations that will be of use in the following.
For any sets $A$ and $B$, we set
\begin{align*}
\check{A}  &  =\{-x,x\in A\},\\
A\ominus B  &  =\{x\in \RR ^2/x+\check{B}\subset A\},\\
A\oplus B  &  =\{x+y,x\in A,y\in B\},
\end{align*}
$A\ominus \check B$ is called the erosion of $A$ by $B$, and $A\oplus \check B$ the dilation of $A$
by $B$.

Since we will need these results in a near future, let us also mention that
$(F,K)\rightarrow F\oplus K$ is a continuous (and thus measurable) function from $\calF\times\calK$ into $\calF$. We also have that $(F,K)\rightarrow F\ominus K$ is
upper semicontinuous from $\calF \times \calK \setminus \{\emptyset\}$ into
$\calF$. These results may be found in \cite{Matheron75}. Using the measurability of
$F\rightarrow\Comp{F}$ from $\calF$ to $\calG$, and $F\oplus
K=\Comp{(\Comp{F}\ominus K)}$, we see that $(G,K)\rightarrow G\ominus K$ and
$(G,K)\rightarrow G\oplus K$ are
measurable maps from $\calG \times \calK$ into $\calG$ (Borelians are defined
on $\calG$ in a way similar to those of $\calF$, see \cite{Matheron75}). In
the sequel, we will compute the probability of events involving such
operations without any further comment.

The dead leaves model is a particular instance of a random
tessellation, obtained through sequential superposition of random objects falling on the plane.
More formally, let $\{(x_i,t_i)\}$ be a homogeneous Poisson point process (P.P.P.) in $\RR^2\times\RR$ with intensity one. Let
$X$ be a random closed set of the plane, and $(X_i)$ be i.i.d. closed sets with the same distribution as $X$, independent of
the Poisson process above. Note that from the assumptions on $\{(x_i,t_i)\}$ and
$(X_i)$, $\Phi:=\{(x_i,t_i,X_i)\}$ is a P.P.P. of
$\RR^3\times\calF$. We first define

\begin{definition}
\label{def:vp}
The random closed set $x_i + X_i$ is called a leaf and
\begin{equation}
  \label{eq:defVi}
V_i=\left(x_i + X_i\right)\setminus\left(\bigcup_{t_j\in(t_i,0)}\left(x_j+\inter{X}_j)\right)\right)  
\end{equation}
is called a visible part.
\end{definition}
Notice that, using the measurability of some standards operations on sets following from \cite{Matheron75}, $\{V_i\}$ is a
collection of random closed sets.
From now on, we assume that $X$ satisfies the following two conditions:
\begin{enumerate}[(C-1)]
\item \label{propdef1}
For all $K \in \K$, $E \nu (X \oplus K) < +\infty$,
\item \label{propdef2}
there exists a disk D with strictly positive radius, such that $E\nu(X \ominus D) >0$.
\end{enumerate}

\begin{proposition}
\label{prop:deadleaves}
$M=\{V_i\}$ is a random tessellation of $\RR^2$.
\end{proposition}

In order to prove this result we will make use of the following two lemmas, the
first of which
will be repeatedly needed in the sequel. 

\begin{lemma}
\label{lemma:reduit}
 Let $K$ be a compact set, $-\infty\leq t_1<t_2<0$ and define 
\begin{eqnarray*}
\Phi_K(t_1,t_2) & = & \left\{(x_i, t_i, X_i) \,:\, t_i\subset (t_1,t_2] \mbox{ and } K \subset x_i+X_i \right\}, \\
\Phi^K(t_1,t_2) & = & \left\{(x_i, t_i, X_i) \,:\, t_i\subset (t_1,t_2] \mbox{ and } K \cap x_i+X_i \neq \emptyset \right\}.
\end{eqnarray*}
Then $\Phi_K(t_1,t_2)$ and $\Phi^K(t_1,t_2)$ are P.P.P.'s of $\RR^3\times\calF$,
respectively of associated measures $\mu_{K,t_1,t_2}$ and $\mu^{K,t_1,t_2}$. Moreover we have
$\mu_{K,t_1,t_2}(\R^3\times\calF)=(t_1-t_2) E\nu(X{\ominus\check{K}})$ and 
$\mu^{K,t_1,t_2}(\R^3\times\calF)=(t_1-t_2) E\nu(X{\oplus\check{K}})$.
\end{lemma}

\paragraph{proof} The result 
follows by a classical result on independent thinning applied to $\Phi$, see
\cite{Pointprocesses2}, and the fact that $(K\subset (x_i+X_i))\Leftrightarrow
(x_i \in X_i\ominus K)$, and $(K \cap (x_i+X_i) \neq \emptyset)\Leftrightarrow
(x_i \in X_i\oplus K)$.
\cqfd

\begin{lemma}
\label{lemma:covering}
If $K\subset\RR^2$, $K$ measurable, is such that $E\nu(X \ominus\check
{K})>0$, then $K$ is covered by some leaf $x_i + X_i$, for some $i$ such
that $t_i\leq 0$, with probability 1. As a consequence, any bounded set is
almost surely covered by a finite number of leaves.
\end{lemma}

\paragraph{proof}
Let us fix $t>0$. Using the same point process $\xx$ as in Lemma
\ref{lemma:reduit}, the probability that none of the leaves $X_i$ with $t<t_i<0$
satisfies $K\subset x_i+ X_i$ is $\exp(-tE\nu(X\ominus \check K))$,
which yields the first assertion. Now let $D$ be a disk such that Condition 
(\textbf{C}-\ref{propdef2}) is satisfied, that is $E
\nu(X \ominus D) >0$. Since any bounded set $K$ is covered by a finite
number of disks with the same radius as $D$, the second assertion follows.
\cqfd

\paragraph{Proof of Proposition \ref{prop:deadleaves}} 
We check that $M$ is almost surely a tessellation. 
From Lemma \ref{lemma:covering}, we know that any compact set is almost
surely covered by a finite number of leaves, which yields
property (i) of Definition \ref{def:tessellation}. Property (ii) of
is directly obtained from~(\ref{eq:defVi}). 
Now, since any compact set $K$ is almost surely covered by a finite number of leaves $V_{i_j}$, $j=1,\ldots,n$, letting
$T=\min(t_{i_j},\,j=1,\ldots,n)$, we see that all leaf $V_j$ with $t_j<T$ has an empty intersection with $K$. 
It follows that the number of leaves intersecting $K$ is less than a Poisson r.v. with intensity
$\mu^{K,T,0}(\R^3\times\calF)=T\,E\nu(X{\oplus\check{K}})$  
(see Lemma~\ref{lemma:reduit}) which is finite by Condition~(\textbf{C}-\ref{propdef1}). 
\cqfd

\begin{remark}
\label{rem:baldleaves}
Recall that in Remark \ref{rem:tessbis}, we pointed out that if the $V_i$'s
are such that $V_i = \overline{\inter{V_i}}$, then the $\inter{V_i}$'s verify
relation (\ref{eq:tess4}). In the case of the dead leaves model, we easily see
that if $X_i = \overline{\inter{X_i}}$, then the $V_i$'s satisfy condition
(\ref{eq:tess4}) (even though in general $V_i \neq \overline{\inter{V_i}}$). We
thus have a way to ensure that the interior of the $V_i$ and the
boundary of our tessellation are disjoint, by ensuring that the $X_i$'s satisfy
(\ref{eq:tess4}). Let us point out that the operation $F\rightarrow
\overline{\inter{F}}$ bears some similarity with the so-called "morphological
opening", defined for $\epsilon>0$ by $F\rightarrow (F\ominus
B(0,\epsilon))\oplus B(0,\epsilon)$,
which is commonly used in image processing to "clean" sets or images. The aim
is the same in our case : we can ensure that our objects do not have
parts with empty interior, such as points or lines.
\end{remark}

We are now in a position to define a dead leaves model associated to a RACS $X$
satisfying (\textbf{C}-\ref{propdef1}) and (\textbf{C}-\ref{propdef2}) as a random tessalation.

\begin{definition}
\label{def:dl}
(The distribution of) $M=\{V_i\}$ is the dead leaves model associated with (the distribution of)
$X$. The set $\partial M$ is defined as in Definition~\ref{def:tessellation} and is called the boundary of $M$. The $V_i$'s
are called the visible parts of M, and connected components of $\R^2\backslash\partial M$ are called the cells of $M$. 
\end{definition}

\begin{remark}
We could also have defined the model through the random closed set $\partial M$. Here we propose a more complete definition
since the tessellation $M$ is generally not uniquely defined by its boundary (see the comments after
Definition~\ref{def:tessellation}).  
\end{remark}

There exist different approaches for defining dead leaves models. 
The choice of
intensity one for the P.P.P. $\{(x_i,t_i)\}$ clearly is not 
restrictive. On the one hand, rescaling the $x_i$'s is equivalent, up to a global
rescaling of the model, to a rescaling of $X$. On the other
hand, any increasing transformation of the $t_i$'s is unimportant as seen from
the definition and, in particular, stopping the superposition of objects at
time $t=0$ is not restrictive. Let us mention that $M$ can also be defined as the stationary distribution of the Markov
process $(M_t)_{t\geq0}$, where $M_t$ is the collection of visible parts associated to the thinned
P.P.P. $\Phi^{\R^2}(-t,0)$. Similarly, for a given compact $K\subset\R^2$, $M\cap K:=\{V_i\cap K\}$,
has the stationary distribution of the Markov chain $(K\cap M_{t_i(K)})_{i\geq0}$, where $(t_i(K))_{i\geq0}$ is defined as the
(almost surely) unique increasing sequence such that $t_0(K)=0$ and for all $i\geq0$, $t_{i+1}(K)$ is the closest $t_j$
strictly below $t_{i}(K)$ such that $(x_j+ X_j)\cap K\neq\emptyset$. Observe that $K\cap
M_{t_{i+1}(K)}$ is obtained by putting the new leaf $x_j+ X_j$ \textit{below} the already fallen leaves. Hence the
limit of $(K\cap M_{t_i(K)})_{i\geq0}$ is hit at the first $i$ such that  $K$ has been completely covered. It turns out that
this $i$ almost surely exists for all $K$ (see Lemma~\ref{lemma:covering}).
The term ``dead leaves model'' originates from a more natural definition which consists in putting each new leaf \textit{above} the
previous ones and then consider the stationary distribution of this Markov chain. It is easily seen that the two definitions
yield the same stationary distribution by a classical Markov Chain argument ( ``coupling from the past'') and that the former
one, that is, the one hitting its limit, may be used for perfect simulation. This elegant argument
was first introduced for the dead leaves model (and more generally for problems of
stochastic geometry) in \cite{Kendall99}, see also \cite{Moller01} and the illustrating web
applet \cite{KendallApplet}. The following simulations of the model have been computed this way, that is, by simulating the chain
$(K\cap M_{t_i(K)})_{i\geq0}$. 
\begin{figure}
\begin{center}
\includegraphics[width=7cm]{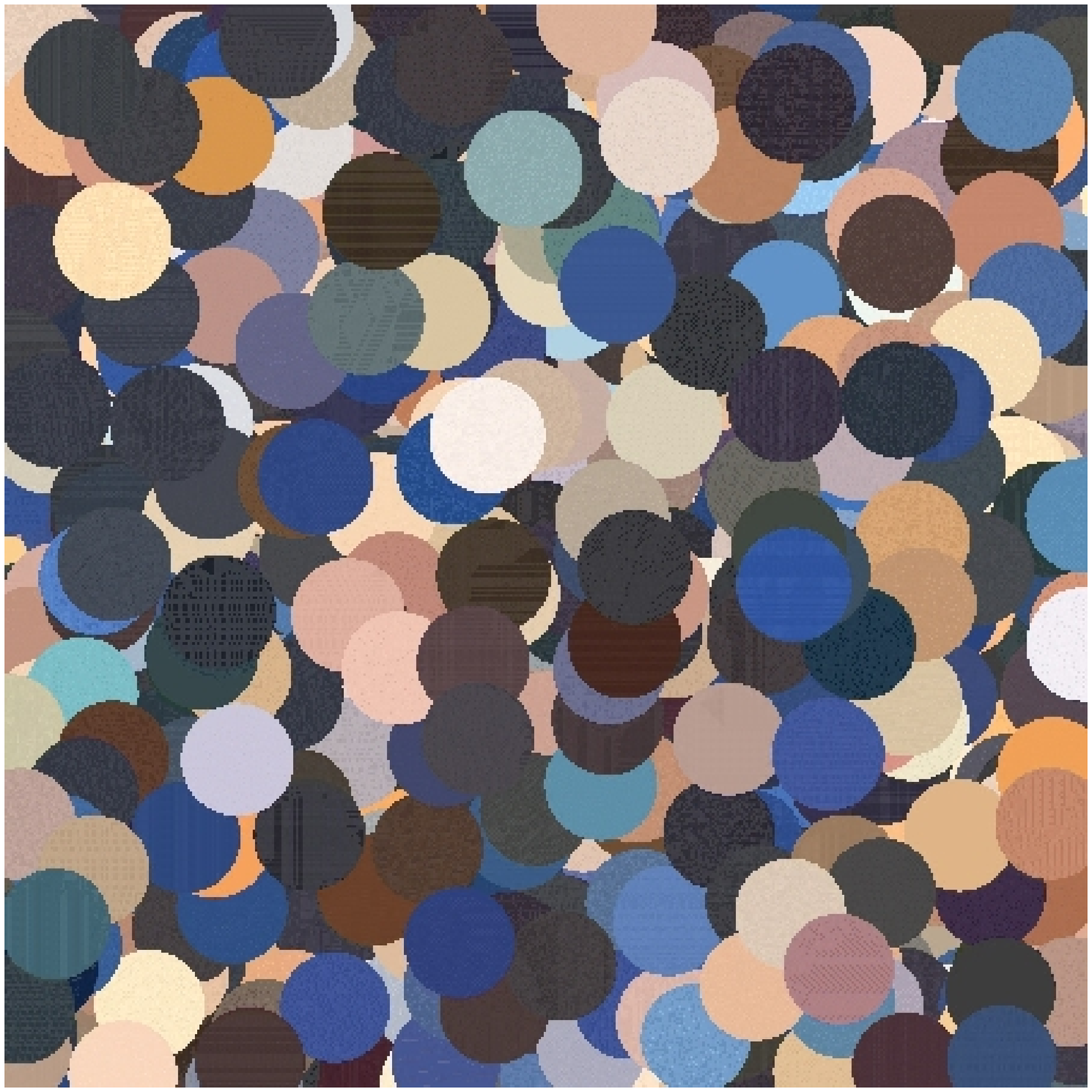} \qquad
\includegraphics[width=7cm]{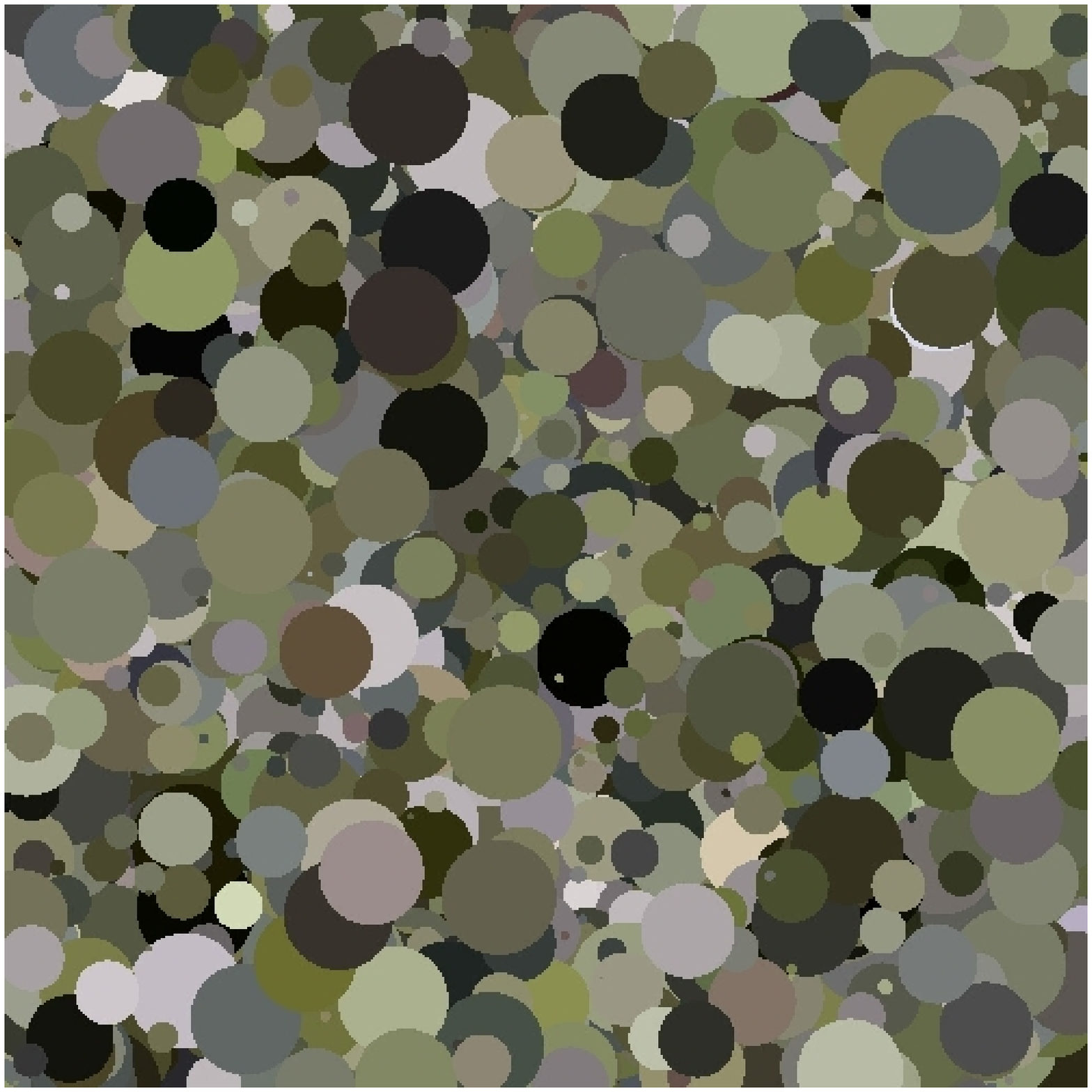}
\caption{Left : simulation of a dead leaves model, where the grain $X_{0}$ is
  a disk with constant radius; Right : simulation of a dead leaves model, where the grain $X_{0}$ is
a disk with a uniformly distributed radius.}
\label{fig:disques}
\end{center}
\end{figure}

\begin{figure}
\begin{center}
\includegraphics[width=7cm]{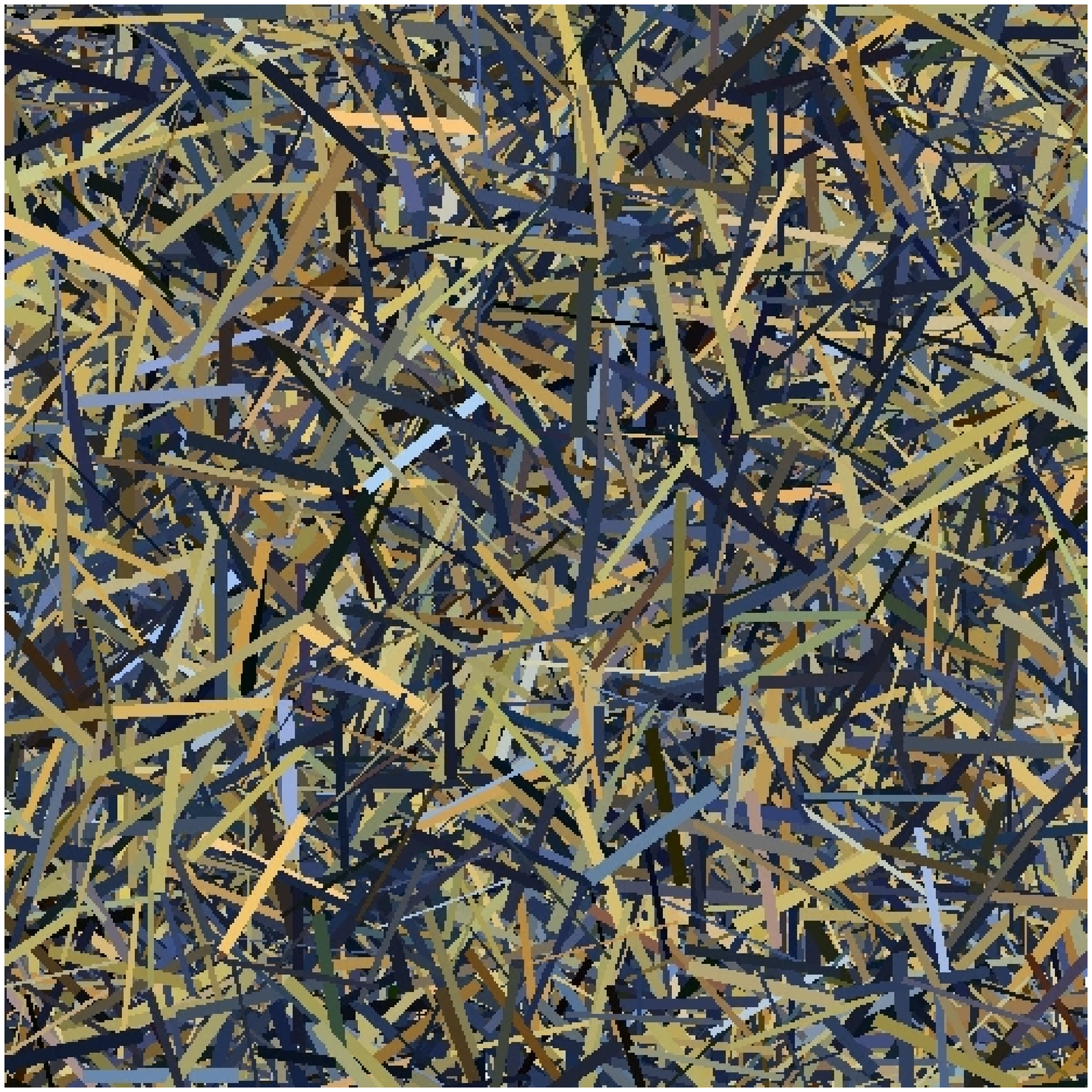} \qquad
\includegraphics[width=7cm]{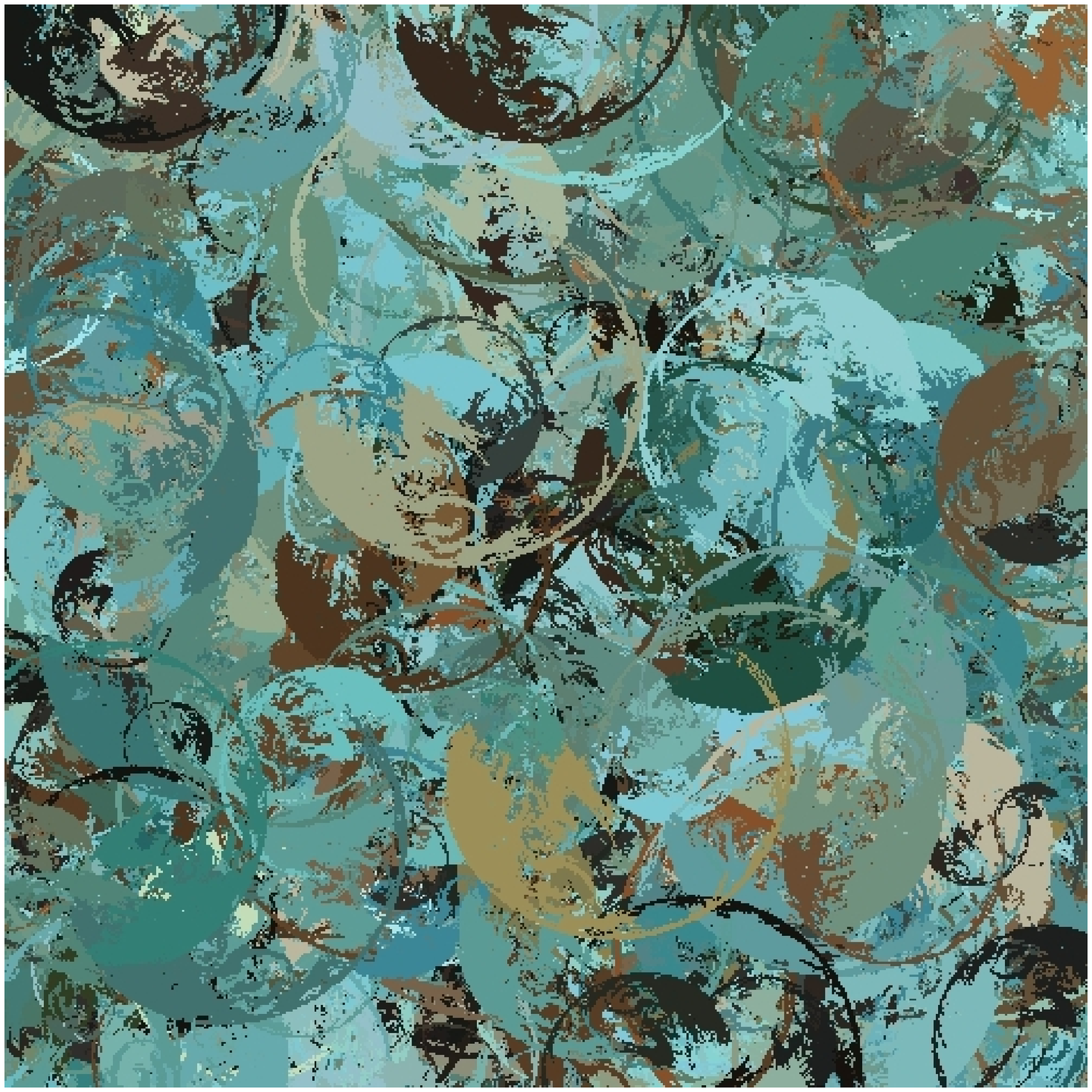}

\caption{simulations of dead leaves models. Left: the grain $X_{0}$ is a
rectangle with a rotation uniformly distributed in $[0,\pi]$. Right: the grain
is more complicated shape, the distribution of its size is uniform.}
\label{fig:recunif}
\end{center}
\end{figure}

In Figures \ref{fig:disques}, and \ref{fig:recunif}, we show realizations of dead leaves models, for different choices of $X$. To
visualize the model each grain is associated a color uniformly drawn in a
table, and the model are ``perfectly'' simulated using the method described in
the previous paragraph. In the first figure, $X$ is a disk of constant radius (left) or
a disk with uniformly distributed radius (right); in the second
one, $X$ is a rectangle whose sides are $r_0 \cos(\theta)$ and
$r_0\sin(\theta)$, with $r_0$ deterministic and $\theta$ uniform between two
bounds, uniformly rotated between 0 and $\pi$ (left) or a more complicated
shape whose size follows a uniform distribution (right).

\subsection{The functional $Q$}

The main practical result from the original paper by Matheron introducing the
dead leaves model ,
\cite{Matheron68}, concerns a functional,
defined on the set of compact sets of the plane, equal to the probability that
a given compact is included in a visible part of the model: for $K \in \K$,
$$
Q(K)=Pr(\exists t_i\leq0, \, K \subset \inter{V}_i).
$$
Considering simple examples of possible $K$ such as bipoints or segments leads
to valuable geometric information on the model.

\begin{proposition}
\label{QKprop}
Let $M$ be the dead leaves model associated with a random closed set $X$, then the probability of $K \in 
\K$ to be included in a visible part of $M$ is
\begin{equation}
\label{QK}
Q(K)=\frac{E\nu(\inter{X}{\ominus\check{K}})}{E\nu(X{\oplus\check{K}})}.
\end{equation}
\end{proposition}

In what follows, we generalize this result by taking interest in the probability that $n$ compact sets are
included in distinct visible parts. In particular this gives us the
probability that any compact consisting of a finite union of connected
components hits the boundary of the dead leaves model. We define 
$$
Q^{(n)}(K_1,\dots,K_n)=Pr(\exists t_{i_1}<\dots<t_{i_n}\mbox{ such that } K_1 \subset \inter{V}_{i_1}, \dots,K_n \subset \inter{V}_{i_n}).
$$

\begin{proposition}
  \label{prop:ncompacts}
Let $M$ be a dead leaves model associated with the random set
$X$, and $\{V_i\}$ its visible parts.   
Let us denote
\begin{align}
\label{eq:ncompacts1}
\begin{split}
F_X^{(n)}(K_1,\dots,K_n)=E\nu(\inter{X}\ominus \check  K_1)
\prod_{j=2}^{n}E\nu\left((\inter{X}\ominus \check  K_j)\cap \Comp{(X\oplus\underline{\check K}_{j-1})}\right),
\end{split}
\end{align}
and 
\begin{equation}
\label{eq:ncompacts2}
G_X^{(n)}(K_1,\dots,K_n)=\prod_{j=1}^{n}E\nu\left(X \oplus \underline{\check K}_{j}\right),
\end{equation}
where, for all $j$,
\begin{equation}
\label{eq:defunderK}  
\underline{K}_j=\bigcup_{k=1}^{j}K_k.
\end{equation}
Then, the probability that $n$ non-empty compact sets
$K_1,\dots,K_n$ be sequentially included in $n$ distinct visible parts is
 
\begin{equation}
  \label{eq:ncompacts}
  Q^{(n)}(K_1,\dots,K_n)= \frac{F^{(n)}_X(K_1,\dots,K_n)}{G^{(n)}_X(K_1,\dots,K_n)}.
\end{equation}
\end{proposition}

\begin{remark}\label{rem:convQn}
We actually see from the formula above that, because $A\ominus B\subset A\oplus C$ for all sets $A,B,C$ such that
$B\subseteq C$, we always have 
$F_X^{(n)}(K_1,\dots,K_n)\leq G_X^{(n)}(K_1,\dots,K_n)$. Note also that 
(\textbf{C}-\ref{propdef2}) implies $E\nu(X)>0$ and thus that $G_X^{(n)}(K_1,\dots,K_n)$
does not vanish for non-empty compact sets.
\end{remark}

\paragraph{proof}

We fix $n$ compacts $K_1,\dots,K_n$. Summing over disjoint events we have that 
$$
Q^{(n)}(K_1,\dots,K_n)=E \left( \sum \indi(t_{i_1}<\dots<t_{i_n})\prod_{j=1}^n\indi(K_{j}
 \subset \inter{V}_{i_j})\right),
$$
where the sum is taken over all $n$-upplets of points in $\Phi$. First note that from the definition of visible parts,
the RHS of this equation may be writen as
$$E\left(\sum \indi(t_{i_1}<\dots<t_{i_n}) \prod_{j=1}^n\indi(K_j \subset (x_{i_j}+\inter{X}_{i_j}))\prod_{t_i>t_{i_j}}\indi(K_j \cap (x_i+X_i) = \emptyset)\right).$$ 
In the simplest case $n=1$, this amounts to say that $Q^{(1)}$ is the probability that there exists a leave $X_i$ such that 
$K_1$ is included in $\inter{X}_i$ and is not hit by subsequent leaves.
We will now apply Campbell theorem to compute this expectation, and therefore need the following notations. Let $\calE := \RR^2 \times \RR \times \F$. We will denote by $z=(x,t,X)$ the points of $\calE$. 
We write $\calN^{(n)}$ ($\calN$ for $n=1$) for the space of point processes on $\calE^n$.
For all $n\geq1$, we then introduce the following point process on $\calE^n$, 
$$
\Phi^{(n)}=\left\{(z_j)_{j=1}^n
  \,:\,  z_1,\,\ldots,z_n\in\Phi \mbox{ disjoints} \right\}.
$$
We consider the following
function, defined for all $\tilde{Z}=(\tilde z_j)_{j=1}^n\in\calE^n$ and for all
$\overline{\Phi}=\{(\overline{x}_{i,j},\overline{t}_{i,j},\overline{X}_{i,j})_{j=1}^n\}\in\calN^{(n)}$ by
\begin{multline}
\label{eq:DefF}
f(\tilde{Z},\overline{\Phi})=\\
\indi(\tilde t_{1}<\dots<\tilde t_{n}) \prod_{j=1}^n\left(\indi(K_{j}
 \subset (\tilde x_j+\inter{\tilde X}_{j})) \prod_{\overline{t}_{i,j}\in(\tilde t_{j},0]}
 \indi(K_{j}\cap (\overline{x}_{i,j}+\overline{X}_{i,j}) =\emptyset)\right). 
\end{multline}
In this equation and below, we use the convention that an indexed variable
overlined with a tilda or a bar denotes a
dummy variable as opposed to those related to $\Phi$.
When applied to $\overline{\Phi}=\Phi^{(n)}$ and $(\tilde z_j)_{j=1}^n=(z_{i_j})_{j=1}^n\in\Phi^{(n)}$, we obtain
$$
Q^{(n)}(K_1,\dots,K_n)=E\left( \sum_{Z\in \Phi^{(n)}} f(Z,\Phi^{(n)})\right).
$$
Applying the refined Campbell theorem (see \cite{Pointprocesses2}), this expectation writes 
$$
Q^{(n)}(K_1,\dots,K_n)=\int_{\calE^n}\int_{\calN^{(n)}}
f(\tilde Z,\Phi^{(n)})\,\tilde P(d \tilde Z)\,\calP_0(\tilde z_1,\dots,\tilde z_n,d\Phi^{(n)}),
$$
where, writing $\tilde Z=((\tilde x_j,\tilde t_j,\tilde X_j))_{j=1}^n$,
$\tilde P (d\tilde Z)=\prod_{j=1}^n (d\tilde x_jd\tilde t_jP(d\tilde X_j))$, $P$ being the probability law of the
random closed set $X$, and where $\calP_0$ is the Palm distribution of the
process $\Phi^{(n)}$. Because $\Phi$ is a
P.P.P. in 
$\calN$, the Slivnyak's Theorem (see \cite{Pointprocesses2}) 
applies, giving
\begin{equation}
\label{eq:slivnyak}
Q^{(n)}(K_1,\dots,K_n)=\int_{\calE^n}\int_{\calN}
f(\tilde{Z}=(\tilde z_j)_{j=1}^n,(\Phi+\delta_{\tilde  z_1}+\dots+\delta_{\tilde z_n})^{(n)})\tilde P(d\tilde Z)\,\calP(d\Phi),
\end{equation}
where $\calP$ is the law of process $\Phi$. For all $\tilde Z=(\tilde{z}_j=(\tilde x_j,\tilde t_j,\tilde X_j))_{j=1}^n$ such that 
$\tilde t_1<\dots<\tilde t_n$, the product in~(\ref{eq:DefF}) may be written as 
\begin{multline}
\label{eq:f2}
f(\tilde Z,(\Phi+\delta_{\tilde z_1}+\dots+\delta_{\tilde z_n})^{(n)})=
\left(\prod_{j=1}^n\indi(K_j\subset \tilde x_j+\inter{\tilde X}_j)\right)\,
\left(\prod_{j=2}^{n} \indi(\underline{K}_{j-1}\cap (\tilde x_j+\tilde X_j)=\emptyset) \right)\\
\left(\prod_{j=1}^{n-1} \prod_{t_i\in(\tilde t_j,\tilde t_{j+1}]}\indi(\underline{K}_j\cap (x_i+X_i) =\emptyset) \right)
\prod_{t_k\in(\tilde t_n,0]} \indi(\underline{K}_n \cap (x_k+X_k) =\emptyset),
\end{multline}
with $\underline{K}_{j}$ as defined in (\ref{eq:defunderK}). To compute the integral in (\ref{eq:slivnyak}), we may first
integrate the term of the second line of (\ref{eq:f2}), first with respect to $\calP$, then to $\indi(\tilde t_{1}<\dots<\tilde t_{n})d\tilde
t_1,\dots,d\tilde t_n$. Now, at fixed $s<t<0$, and for $K$ compact, we have that
$$
\calP(K\cap(x_i+X_i)=\emptyset \mbox{ for all }t_i\in(s,t])=\exp\left((s-t) E \nu (X\oplus \check K)\right),
$$
which follows from  Lemma~\ref{lemma:reduit}. 
Then using a change of variable
$u_j=\tilde t_{j}-\tilde t_{j+1}$, for $j=1...n-1$, we obtain 
\begin{multline}
Q^{(n)}(K_1,\dots,K_n)= 
\prod_{j=1}^{n}
E\nu\left(X \oplus \underline{\check K}_{j}\right)^{-1}\\ 
\int_{(\R^2\times \calF)^n} 
\left(\prod_{j=1}^n\indi(K_j\subset  \tilde x_j+\inter{\tilde X}_j)\right) 
\left(\prod_{j=2}^{n} \indi(\underline{K}_{j-1}\cap (\tilde x_j+\tilde X_j)=\emptyset) \right) 
\prod_{j=1}^n(d\tilde x_jP(d\tilde X_j)). 
\end{multline}
The first term of the RHS of the previous equation is $(G_X^{(n)})^{-1}$, and the term of the second line writes
\begin{equation*}
\prod_{j=1}^{n} \left(\int_{\R^2\times \calF} \indi(K_j\subset\tilde x+\inter{\tilde X}) \indi( \underline{K}_{j-1}\cap
  (\tilde x+\tilde X)=\emptyset)dxP(d\tilde{X})\right),
\end{equation*}
with the convention $\underline{K}_{0}=\emptyset$. Now, for two compact sets $A$ and $B$, we have
$$
\int \indi(A\subset (x+\inter{X}))\indi(B\cap (x+X)=\emptyset) dxP(dX) = 
E_X \nu ((\inter{\check X} \ominus A) \cap (\check X \oplus B)^c ) ,
$$
which, along with the last equations, yields (\ref{eq:ncompacts}). \cqfd

In case $n=1$, we get the original result of Matheron, Proposition
\ref{QKprop}, and the case $n=2$ was treated in \cite{Jeulin}. Note also that we can compute the probability for $n$
connected compact sets $K_1,\dots,K_n$ to avoid the
boundary of the dead leaves model $M$. In case $n=2$ for instance, this is
$$Pr((K_1\cup K_2) \cap \partial
M = \emptyset)=Q^{(2)}(K_1,K_2)+Q^{(2)}(K_2,K_1)+Q^{(1)}(K_1\cup K_2).$$
More generally, we show that the knowledge of $Q^{(n)}$ for all $n$
uniquely determines the probability for any compact set $K$ to avoid the
boundary of $M$, which characterizes its distribution in $(\calF,\calB_\calF)$.

\begin{proposition}
 Let $M$ be a dead leaves model, then the distribution of the random closed set $\partial M$ is
 uniquely determined by the functionals $Q^{(n)}$, $n \in \NN$. 
\end{proposition}

\paragraph{proof}
Thanks to Proposition \ref{prop:CaracRACS}, it is enough to prove that the
functional $Q(K)=Pr(K\cap \partial M = \emptyset)$, defined for compact sets
$K$, is uniquely determined by the $Q^{(n)}$. Let $K \in \K$, let $r_n>0$
be a
sequence converging to 0,
and for each $n$, let $\{x_i^{(n)}\}_{i=1,\dots,N_n}$ be finite sequences in $K$ such
that $K\subset C_n=\cup_i D(x_i^n,r_n)$, where $D(x,r)$ is the (closed) disk centered at
$x$ with radius $r$. Note that since each $C_n$ is a finite union of connected
compact sets, the knowledge of the $Q^{(i)}$, $i\in \NN$, uniquely determines
$Q(C_n)$. Now since $C_n \downarrow K$, we have that $\F^{C_n} \uparrow \F^K$,
and thus that $Q(C_n)\uparrow Q(K)$. \cqfd

\section{A dead leaves model with scaling properties}
\label{scaling}

In this chapter, we take interest in a dead leaves model that has a specific
object distribution. $X$, the grain of the model, is equal to $RY$,
the homothetic of a random compact set $Y$, where $R$ is a random variable
with density $f$, independent of $Y$. We will consider in details the case 
$f(r)=Cr^{-\alpha}\indi(r_0\leq r\leq r_1)$, with $\alpha>1$, that is when the size distribution of the objects
satisfies some scaling properties within in a given range. As explained in the
introduction, this choice for $f$ is motivated by natural images modeling. Clearly we cannot take $r_0=0$ for $f$ to be a density. 
Now, from a modeling point of view, taking $r_0>0$ is not satisfying as well because this would be inadequate for capturing the presence of small objects in natural images.
From a theoretical point of view, this reduces the model to only very simple smoothness classes (namely, piece-wise
constant images). From a practical point of view, it means that there exists a minimal size for the
objects in the image. It is not clear at all what physical meaning to give to
this minimum object size, and how to deal with this supplementary parameter of
the model. It is also unclear how to relate this minimum size to the
resolution of an image, which we may assume to be obtained from
the model through filtering and subsampling. Moreover, this contradicts
empirical experiments (see \cite{Gousseau01}) which conclude to the presence
of small "objects" up to the smallest observable scales in digital images.
Therefore it is worthwhile to wonder about the limit of the model as $r_0$ tends to zero. 
The parameter $r_1$ is not crucial for modeling smoothness properties
because it does not influence the small scales behavior, except perhaps when $r_1=\infty$ or $r_1$ tends to $\infty$ in
which case the model may degenerate. 
These cases will also be discussed. In the \textit{good cases} we will consider both the convergence of
the model boundary set and the finite-dimensional convergence of a colored version of the model, giving raise
to quite different kinds of limits. 
 
\subsection{Definitions}\label{sec:definitions}

Let $Y \subset \R^2$ be a random compact set, and, for any $0<r_0< r_1$, define the probability density function 
\begin{equation}
  \label{eq:defDensGeo}
f(r_0,r_1,r)=\eta(r_0,r_1)\indi(r_0\leq r\leq r_1)r^{-\alpha},
\end{equation}
with $\eta(r_0,r_1)=(1-\alpha)^{-1}(r_0^{1-\alpha}-r_1^{1-\alpha})$, $f=0$ otherwise.
We consider $M(r_0,r_1)$ the dead leaves model associated with the random closed
set $X=R Y$, where $R$ is a random variable with density $f(r_0,r_1,\cdot)$. For convenience, our notations do not refer to
the scaling parameter $\alpha$. However, it must be kept in mind that these definitions highly depend on this parameter.
In order for this model to make sense, we want to keep $Y$ within fixed
proportions. All along the paper it is  assumed that

\begin{hyp}{propy} 
there exist $a_2>a_1>0$ such that almost surely $D(a_1) \subset Y \subset D(a_2)$,
\end{hyp}
where $D(a)$ is the disk of radius $a$ centered at the origin.

It is straightforward to check that $X$ defined this way satisfies the requirements
for $M(r_0,r_1)$ to be a dead leaves model, (C-\ref{propdef1}) and
(C-\ref{propdef2}) for all $0<r_0<r_1<\infty$ and, if $\alpha>3$, also for $r_1=\infty$. 

Let $p(r_0,r_1,\xx)$ denote the probability that the origin and $\xx\in\R^2$ are in the same visible part of the 
dead leaves model $M(r_0,r_1)$. In the next section we focus on the limit behavior of $p(r_0,r_1,\cdot)$ as $r_0$ and $r_1$
goes to 0 and $\infty$ respectively. This will tell us under which condition for $\alpha$ a non-degenerate limit process
can be obtained. According to Formula (\ref{QK}), for all $\xx\in \R^2$, 
\begin{equation}
\label{eq:pxx}
p(r_0,r_1,\xx) = \frac{E \nu\left(\inter X \ominus \{0,\xx\}\right)}{E \nu\left(X \oplus 
\{0,{\xx}\}\right)} = \frac{E \nu \left(\inter X \cap ({\xx}+ \inter X)\right)}{E \nu \left(X \cup ({\xx}+
X)\right)},
\end{equation}
where $\nu$ is the 2-dimensional Lebesgue measure, and $E$ the expectation
with respect to the law of $X$. Writing $E_Y$ for the expectation with respect
to the law of $Y$, Fubini's Theorem and then a geometric argument along with~(\ref{eq:defDensGeo}) 
give
\begin{align}\nonumber
 E \nu\left(X \cap \xx+X\right) & = \int_{r_0}^{r_1}E_Y \nu\left(u\inter Y \cap (\xx+u\inter Y)\right) \, f(r_0,r_1,u)\,du\\
\label{eq:homog}
&= \eta(r_0,r_1)\,\int_{r_0}^{r_1} \tilde\gamma\left(\frac{\xx}{u}\right)\,u^{2-\alpha}\,du, 
\end{align}
where $\tilde\gamma$ denotes the geometric covariogram of
$Y$, see \cite{Matheron75}, that is, for all $\yy\in\R^2$,
$$
\tilde\gamma(\yy):=E_Y \nu\left(\inter Y \cap (\yy+\inter Y)\right).
$$
Let us now consider some basic assumptions and results concerning the mapping $\tilde\gamma$.
From \refhyp{propy}, we obtain
\begin{lemma}\label{lem:gammatilde} 
The mapping $\yy\mapsto\tilde\gamma(\yy)$ is continuous over $\R^2$,
$\tilde\gamma(0)\geq\pi a_1^2$ and, for all $\yy\in\R^2$, 
$0\leq\tilde \gamma(\yy)\leq \tilde{\gamma}(0)$, where $|\cdot|$ is the
Euclidean norm. Moreover,
\begin{equation}
  \label{eq:gammatildeinfty}
|\yy|\geq2a_2\Rightarrow \tilde \gamma(\yy)=0.
\end{equation}
\end{lemma}
\paragraph{proof}
The bounds on $\tilde\gamma$ are immediate. We now notice that $\yy\mapsto\nu\left(\inter Y \cap (\yy+\inter Y)\right)$ is the
convolution of the indicator function on $\inter Y$ with itself. Since $Y$ is bounded, this indicator
function is square integrable with respect to $\nu$ and the convolution is continuous. By the dominated convergence theorem,
the continuity is preserved after taking the expectation.
\cqfd 

The two following lemmas show that a stronger regularity of the
geometric covariogram at the origin may be obtained by imposing some geometrical properties on $Y$. First it is known,
see \cite[Proposition 4-3-1]{Matheron75}, that   
\begin{lemma}
\label{gammaconvex}
Let $\bu$ be a unitary vector. Let $C$ be a closed convex set such that $\inter C \neq 
\emptyset$, then $x\mapsto\nu(C\cap C_{x\bu})$
admits a derivative to the right at $x=0$ equal to $\nu_1(P_{\bu}(C))$,
where $\nu_1$ is the one dimensional Lebesgue measure, and $P_\bu$ the
orthogonal projection in the direction of $\bu$. 
\end{lemma}

This lemma can for instance be generalized to finite unions of convex sets. 
Simple criteria on the random set $Y$ can thus be deduced in ordered to impose
that 

\begin{hyp}{gamma2} 
  for every unitary vector  $\bu$, the function $\tilde \gamma (x\bu)$ is
  differentiable to the right at $x=0$ and the derivative is bounded away from the origin and infinity independently of
  $\bu$.  
\end{hyp}

This condition will be used as an assumption in some of the results below. However, in most cases we will relax this
assumption into the following one.

\begin{hyp}{gamma1} 
  For any $\delta>0$, we have $\tilde \gamma(\xx)=\tilde{\gamma}(0)+o(|\xx|^{1-\delta})$ when
  $\xx \rightarrow 0$. 
\end{hyp}

For obtaining such conditions, we have the following bound.

\begin{lemma}
  \label{gammarectif}
Let $C$ be a compact set. Then for all $\xx\in\R^2$,
$$
0\leq\nu( \inter C)-\nu(\inter C\cap(\xx+\inter C))\leq\nu\left(\bigcup_{\bu\in\partial{C}}(\bu+[0,\xx])\right),
$$
where $[0,\xx]$ denotes the segment of points between $\xx$ and the origin, that is $\{\alpha\xx,\,\alpha\in[0,1]\}$. 
\end{lemma}
\paragraph{proof} 
Let $\xx\in\R^2$. We may write 
$$
0\leq\nu(\inter C)-\nu(\inter C\cap(\xx+\inter C))=\nu(\inter C\backslash(\xx+\inter C))
$$
Let $\yy\in \inter C\backslash(\xx+\inter C)$, that is $\yy\in \inter C$ and $\yy-\xx\notin \inter C$. Thus the line segment
$[\yy,\yy-\xx]$ intersects $\partial C$. The result follows.
\cqfd

In particular we see that if $\partial C$ is a rectifiable 
curve with length $L$, then
$0\leq\nu(\inter C)-\nu(\inter C\cap(\xx+\inter C))\leq L|\xx|$. More generally, if $C$ is such that its topological boundary, $\partial C$     
has an upper box-counting dimension less than or equal to one, then
(see \cite[Proposition 3.2]{Falconer90}), for all $\delta>0$,  
$$
\nu\left(\bigcup_{\bu\in\partial C}D(\bu,a)\right)=o(a^{1-\delta}),
$$
where $D(\xx,a)$ denotes the disk of radius $a$ centered at $\xx$. By Lemma~\ref{gammarectif}, the same bound applies to 
$\nu(\inter C)-\nu(\inter C\cap(\xx+\inter C))$ with $a=|\xx|$ in a case where the boundary of $C$
is not necessarily a parameterized
rectifiable curve. Hence Assumption~\refhyp{gamma1} is easy to verify for a large
class of random closed sets $Y$ if one has some minimal property on the measure of $\partial Y$. This is summarized in the
following corollary. 
\begin{corollary}\label{cor:assumptions}
Assume that for any $\delta>0$, we have
$$
E_Y\nu\left(\bigcup_{\bu\in\partial{Y}}(\bu+[0,\xx])\right)=o(|\xx|^{1-\delta})
$$
when $\xx \rightarrow 0$. Then~\refhyp{gamma1} holds true.
\end{corollary}

However, in order to illustrate that some smoothness is needed on the boundary of $Y$
to insure \refhyp{gamma1}, even for a deterministic $Y$. Consider the following example. Let $h:[0,1]\to(0,1)$ be a
continuous mapping and let us define
$$
C:=\left\{(x,y):0\leq x\leq 1,\,h(x)-1\leq y\leq h(x)\right\}.
$$
Then we have, for all $u>0$,
$$ 
\nu(\inter C)-\nu(\inter C\cap((u,0)+\inter C))=\nu(\inter C\backslash((u,0)+\inter C))
=u+\int_u^1|h(x)-h(x-u)|\,dx. 
$$
Hence in this case, taking say $Y=C$ non-random (it could be made random by taking $h$ random), \refhyp{gamma1}
would imply that, for all positive $\delta$, $\int_u^1|h(x)-h(x-u)|\,dx=o(u^{1-\delta})$ as $u\downarrow0$. This
exactly means that $h\in \bigcap_{\delta<1} B_1^{\delta,\infty}([0,1])$ (see
e.g. \cite{devore91} or \cite{meyer90} and Section~\ref{prop:BesovEsp}, where
the Besov spaces $B_p^{s,q}$ are defined) and continuous bounded mappings may be
easily found out of $B_1^{\delta,\infty}([0,1])$ for some $\delta<1$. In such cases 
\refhyp{gamma1} cannot be satisfied. This example shows that this assumption amounts to some smoothness assumption
on the boundary of $Y$.

Another way to interpret Assumption~\refhyp{gamma1} is to say that in the limit model we shall consider, the
smoothness is mainly driven by objects sizes distribution and not by the irregularity that may occur in
their boundaries.

\subsection{Basic convergence results}\label{sec:basic-conv-results}

In order to simplify the presentation of this section, we temporarily assume that the distribution of $Y$ is
isotropic. However, all the convergence results of this  
section remains true in the case of nonisotropic $Y$, the various quantities under study
then depending on a directional parameter. In the isotropic case, we let $\tilde\gamma$ and $p(r_0,r_1,\cdot)$
be mappings of the real variable $x=|\xx|$. The general case may be obtained in the same way by letting these mapping also
depend on the angle coordinate of $\xx$. 

From~(\ref{eq:homog}) we may investigate under which conditions $p(r_0,r_1,\cdot)$ has a
non-degenerate limit when one wants to push the model towards the values of $r_0$ and $r_1$ 
which are not allowed, that is, for all $\alpha>1$, $r_0=0$ because $f(0,r_1,\cdot)$ is not well defined
and, for all $1<\alpha\leq3$, $r_1=\infty$, in which case, condition~(C-\ref{propdef1}) does not hold. 

The cases for which the obtained limits degenerate are summarized in the following Proposition.
\begin{proposition}
\label{deg_limits_synthese}
For all $x>0$,
\begin{enumerate}[{\rm (i)}]
\item \label{eq:deg1} if $1<\alpha<3$, $\displaystyle\lim_{r_1\to\infty}\inf_{0<r_0<r_1}p(r_0,r_1,x)=1$,
\item \label{eq:deg2} if $\alpha>3$, $\displaystyle\lim_{r_0\to0}\sup_{0<r_0<r_1}p(r_0,r_1,x)=0$,
\item \label{eq:deg3} if $\alpha=3$,
  $\displaystyle\lim_{\stackrel{r_0\to0}{r_1\to\infty}}\left[p(r_0,r_1,x)\left(1-2\frac{\log(r_0)}{\log(r_1)}\right)\right]=1$. 
\end{enumerate}
\end{proposition}

\paragraph{proof}
Using~(\ref{eq:homog}) and since 
$\nu \left(X \cup ({\xx}+X)\right)=2\nu\left(X\right) - \nu\left(X \cap \xx+X\right)$, Eq.~(\ref{eq:pxx}) reads
\begin{equation}
\label{eq:Pgamma}
p(r_0,r_1,\xx) = =\frac{
\int_{r_0}^{r_1} \left[\tilde{\gamma}(0) + (\tilde{\gamma}(\xx/u)-\tilde{\gamma}(0))\right] u^{2-\alpha}\,du}
{\int_{r_0}^{r_1} \left[\tilde{\gamma}(0) - (\tilde{\gamma}(\xx/u)-\tilde{\gamma}(0))\right] u^{2-\alpha}\,du}. 
\end{equation}
We now derive the asymptotic behavior of the integrals in this formula depending on the value of $\alpha>1$.

Take $1<\alpha<3$. Let $x>0$. 
From~(\ref{eq:gammatildeinfty}) and since
$\tilde{\gamma}(0)>0$, for any $\epsilon>0$, there exists $u_0$ such that for all $u\in[u_0,\infty)$, 
$$
|\tilde{\gamma}(x/u)-\tilde{\gamma}(0)|\leq \epsilon \tilde{\gamma}(0).
$$
Since the integral below diverges as $r_1\to\infty$ and is bounded as $r_0\to0$, for $r_1$ sufficiently large and for
all $r_0\leq r_1$,  
$$
\int_{r_0}^{r_1}u^{2-\alpha}\,du\leq (1+\epsilon) \int_{u_0\vee r_0}^{r_1}u^{2-\alpha}\,du.
$$
From the two last equations and from~(\ref{eq:Pgamma}), we get, for all $r_1$ sufficiently large and for all $0<r_0<r_1$, 
$$
1\geq p(r_0,r_1,x) \geq \frac{\tilde{\gamma}(0)(1-\epsilon)\int_{u_0\vee r_0}^{r_1}u^{2-\alpha}\,du}
{\tilde{\gamma}(0)(1+\epsilon)^2\int_{u_0\vee r_0}^{r_1}u^{2-\alpha}\,du}=\frac{1-\epsilon}{(1+\epsilon)^2}. 
$$
Hence~(\ref{eq:deg1}) by letting $\epsilon$ decrease to zero.

Take now $\alpha>3$. From~(\ref{eq:gammatildeinfty}) and~(\ref{eq:Pgamma}),
we have, for all $r_0<r_1\leq x/(2a_2)$,  $p(r_0,r_1,x)=0$ and, for all $r_0\leq x/(2a_2)< r_1$,
$$
\int_{r_0}^{r_1} \tilde{\gamma}(x/u) u^{2-\alpha}\,du=\int_{x/(2a_2)}^{r_1}
\left[\tilde{\gamma}(x/u)\right] u^{2-\alpha}\,du.
$$
Since $\tilde\gamma$ is bounded by $\pi a_2$ from above, we get, for all $r_0\leq x/(2a_2)$ and for all $r_1>r_0$, 
$$
p(r_0,r_1,x)\leq \frac{\pi a_2 \int_{x/(2a_2)}^{\infty}u^{2-\alpha}\,du}
{2\tilde\gamma(0)\int_{r_0}^{x/(2a_2)}u^{2-\alpha}\,du-\pi a_2 \int_{x/(2a_2)}^{\infty}u^{2-\alpha}\,du}
$$
which does not depend on $r_1$ and tend to zero as $r_0\to0$. This gives~(\ref{eq:deg2}).

We now conclude with the case $\alpha=3$. From~(\ref{eq:gammatildeinfty}) and the continuity of $\tilde\gamma$, the
numerator of the RHS of~(\ref{eq:Pgamma}) behaves as $\tilde\gamma(0)\log(r_1)$ when $r_0$ and $r_1$ respectively tend
to 0 and $\infty$. For the same reasons, the denominator behaves as $\tilde\gamma(0)(\log(r_1)-2\log(r_0))$. We finally
obtain~(\ref{eq:deg3}). 
\cqfd

It is worth commenting on the consequence of these simple convergence results. In case~(\ref{eq:deg1}), the result
says that, as $r_1$ tends to infinity, however $r_0<r_1$ may behave, any two points end up in the same visible part. The
big objects predominate at the limit so that the image is simply covered by one
single set $X$ with probability one. In case~(\ref{eq:deg2}), the result is the exact opposite. As $r_0$ tend to zero,
however $r_1\in(r_0,\infty]$ may behave, any two points never belong to the same object. The small objects predominate so
that every visible part reduces to a point and the limit distribution of the image is white noise. See Figure \ref{fig:degenerate}
for an illustration of these cases. Finally, in the
case~(\ref{eq:deg3}), the limit depend on the behavior of $\log(r_0)/\log(r_1)$.
Convergence to 1 or 0 as in cases~(\ref{eq:deg1}) and~(\ref{eq:deg2}) are observed if only one of the
limit $r_0\to0$ or $r_1\to\infty$ is taken. If, for instance we take $r_0=r_1^{-s}$ for a fixed $s$, and let $r_1$ tend to
$\infty$, we obtain a limit which depend on $s$ but does not depend on $x$. 

We will generally avoid these cases in the sequel as they obviously give degenerate limits. In contrast, when $r_0\to0$
for $1\leq \alpha<3$, the limit of $p$ is easily obtained from~(\ref{eq:Pgamma}) giving raise to the continuous
prolongation 
\begin{equation}
  \label{eq:prolong_p_en0}
  p(0,r_1,x)  =
=\frac{
\int_{0}^{r_1} \tilde{\gamma}(x/u) u^{2-\alpha}\,du}
{\int_{0}^{r_1} \left[2\tilde{\gamma}(0) - \tilde{\gamma}(x/u)\right] u^{2-\alpha}\,du}.
\end{equation}
These integrals trivially converge since $\tilde{\gamma}$ is bounded and $\alpha<3$. In the following result, by
assuming sufficient smoothness on the boundary of $Y$, we exhibit
asymptotical behaviors in which the geometry of the model only appear in multiplicative constants while the qualitative
behavior is a power law of the form $x^{3-\alpha}$. 

\begin{proposition}\label{prop:asympP}
We have the following asymptotic expansions.
\begin{enumerate}[{\rm (i)}]
\item\label{item:asym1} For all $\alpha>3$, $p(r_0,\infty,x)=g(\alpha)\,(x/r_0)^{3-\alpha} (1+o(1))$ as
  $x/r_0\to\infty$, where
$$
g(\alpha):=\frac{\alpha-3}{2\tilde\gamma(0)}\int_0^\infty\tilde\gamma(1/v)\,v^{2-\alpha}\,dv<\infty
 \mbox{ for all }\alpha>3.
$$
\item\label{item:asym2} Under~\refhyp{gamma1}, for all $2<\alpha<3$,
  $1-p(0,r_1,x)=g(\alpha)\,(x/r_1)^{3-\alpha}(1+o(1))$ as $x/r_1\to0$, where
$$
g(\alpha):=\frac{2(3-\alpha)}{\tilde{\gamma}(0)}\int_0^\infty(\tilde\gamma(0)-\tilde\gamma(1/v))\,v^{2-\alpha}\,dv<\infty 
\mbox{ for all }2<\alpha<3
$$
\item\label{item:asym3} Under~\refhyp{gamma2},  for $\alpha=2$, 
$\displaystyle 1-p(0,r_1,x)=\frac{2\dot{\tilde \gamma}(0)}{\tilde{\gamma}(0)} 
(x/r_1)\log(x/r_1)(1+o(1))$ as $x/r_1\to0$,
\item\label{item:asym4} Under~\refhyp{gamma2}, for all $1<\alpha<2$, 
$\displaystyle 1-p(0,r_1,x)=g(\alpha)(x/r_1)(1+o(1))$ as
$x/r_1\to0$, where
$$
g(\alpha):=\frac{2(3-\alpha)\dot{\tilde \gamma}(0)}{(\alpha-2)\tilde{\gamma}(0)}<\infty
\mbox{ for all }1<\alpha<2.
$$
\end{enumerate} 
In the formula above, under~\refhyp{gamma2} and having assumed that $Y$ has an isotropic distribution, $\dot{\tilde
  \gamma}(0)$ denotes the right-sided derivative of $x\mapsto\tilde\gamma(x)$, $x\in[0,\infty)$, at the
  origin. 
\end{proposition}

\paragraph{proof}
Take $\alpha>3$. An obvious change of variable gives 
$$
\int_{r_0}^\infty \tilde\gamma(x/u)u^{2-\alpha}\,du = 
x^{3-\alpha}\int_{r_0/x}^\infty \tilde\gamma(1/v)v^{2-\alpha}\,dv.
$$
Using~(\ref{eq:gammatildeinfty}) and~(\ref{eq:Pgamma}), letting $x$ decrease to zero gives~(\ref{item:asym1}).
We now take $\alpha<3$. We similarly have
$$
\int_0^{r_1}(\tilde\gamma(x/u)-\tilde\gamma(0))u^{2-\alpha}\,du
=x^{3-\alpha}\int_0^{r_1/x}(\tilde\gamma(1/v)-\tilde\gamma(0))v^{2-\alpha}\,dv.
$$
From~(\ref{eq:Pgamma}) and~(\ref{eq:prolong_p_en0}) and standard computations, we
obtain~(\ref{item:asym2})-(\ref{item:asym4}). 
\cqfd

It is worth noting that empirical experiments have shown that natural images exhibit such power law behaviors, with
$\alpha$ having values usually in the range $[2.5,3]$ (see \cite{Ruderman96},
and also various study on the power spectrum of natural images). Hence in the following we will focus on
cases~(\ref{item:asym1}) and~(\ref{item:asym2}). 
One way to interpret these two cases is the following. The behavior at large scale of $p(r_0,\infty,x)$ with $r_0>0$ and
$\alpha>3$, is qualitatively similar (namely a power law) to the behavior at small scales of $p(0,r_1,x)$ with
$r_1>0$ and $2<\alpha<3$. Hence the limit of $M(r_0,r_1)$ as $r_0\to0$ for a fixed $r_1>0$ seems the natural
way to extrapolate the model $M(r_0,r_1)$ from $r_0>0$, $r_1=\infty$ and $\alpha>3$ to $r_0=0$, $r_1>0$ and
$\alpha<3$. Based on this observation, we dedicate the next sections to defining a limit model $M(0,r_1)$ for $2<\alpha<3$.
It turns out this limit is degenerate from a random set point of view (see Section~\ref{racs}) but it is not from a
colored model point of view (see Section~\ref{colored}).

\subsection{Limit of the boundary set} \label{sec:limit-boundary-set}
In this section we assume that 
$1<\alpha\leq3$ and $r_1>0$ or that $\alpha>3$ and $r_1\in(0,\infty]$. 
We have seen in Section~\ref{sec:basic-conv-results} that, 
although the model
$M(0,r_1)$ is not correctly defined, the probability $p(r_0,r_1,x)$ can be continuously prolonged at $r_0=0$ for all $x$.
In this section, we investigate how the model $M(r_0,r_1)$ \textit{converges} under the same conditions.
Here we take the classical point of view of random closed sets and we consider the distribution of the (random) boundary
set. 
In the case $\alpha>3$, $p(0,r_1,x)=1$ (see Proposition~\ref{deg_limits_synthese}) so that we expect a degenerate limit
of the dead leaves model. 
It turns out that the limit is also degenerate in the case $1<\alpha\leq3$ as shown by
the following proposition.

\begin{proposition}\label{prop:BoundaryLimit}
Let $\partial M(r_0,r_1)$ denote the
random set consisting of the union of the boundaries of visible parts in the model $M(r_0,r_1)$. Then 
$$\partial M \underset{r_0\rightarrow 0}\rightarrow \RR^2,$$
the convergence being the weak convergence of random closed sets.  
\end{proposition}

This convergence result follows easily from the next proposition which investigates the presence of constant areas in a  
possible limit process as $r_0\rightarrow0$ at fixed $r_1$.

\begin{proposition}
\label{InclDisque1}
Let  $Q(r_0,r_1,r)$ denote the probability for a disk of radius $r$
to be included in a visible part of $M(r_0,r_1)$. Then, for any $r>0$,
$\displaystyle\lim_{r_0\rightarrow 0} Q(r_0,r_1,r)=0.$
\end{proposition}

\paragraph{proof}
According to formula (\ref{QK}) and then to \refhyp{propy}, we have
$$
Q(r_0,r_1,r)=\frac{E\nu(\inter X\ominus  D(r))}{E\nu(X\oplus D(r))}
\leq\frac{E\nu(R D(a_2)\ominus D(r))}{E\nu(R D(a_1)\oplus D(r))}
=\frac{\int_{a_2^{-1}r}^{r_1}\pi (u a_2-r)^2\,u^{-\alpha}\,du}{\int_{r_0}^{r_1}\pi (u a_1+r)^2\,u^{-\alpha}\,du}. 
$$
The limit is now obvious.\cqfd

\paragraph{proof of Proposition~\ref{prop:BoundaryLimit}}
Let $P(r_0,r_1,\cdot)$ denote the probability law of $\partial M(r_0,r_1)$ in the probability space $(\F,B_{\F})$ (see
Section~\ref{racs}). We
recall that a sequence $P_n$ weakly converges to $P$ 
in $(\F,B_{\F})$ if for all $E \in B_{\F}$ such that $P(E)=P(\inter{E})$, $P_n(E)$ converges to $P(E)$ (see \cite{Billingsley68}). Moreover, in the case of the probability space $(\F,B_{\F})$, this amounts to check that for all $K\in\K$ such that
$P(\F_K)=P(\F_{\inter{K}})$, $P_n(\F_K)$ converges to $P(\F_K)$ (see
\cite{Lyashenko}, \cite{MolchanovLNM}). 
Here the limit distribution $P$ associated with the deterministic set $\R^2$ satisfies $P(\F_K)=1$ for all compact set
$K\neq\emptyset$ and $P(\F_\emptyset)=0$. Take a compact set $K$ such that $\inter{K}\neq \emptyset$.  
There exist a disk with positive radius $r$ included in $K$ so that
$$
P(r_0,r_1,\F_K)\geq P(r_0,r_1,\F_{D(r)})=1-Q(r_0,r_1,r).
$$
The result then follows from Proposition~\ref{InclDisque1}.\cqfd

\subsection{The colored dead leaves process and its limit}

For modeling the effect of the occlusions on the smoothness in natural images in a non-trivial way, we consider the limit of
the model $M(r_0,r_1)$ as $r_0\to0$. From the point of view of  random closed sets, the limit degenerates
as we have seen in Proposition~\ref{prop:BoundaryLimit}. The
possible limit set is not well described by the boundary limit, and we now take interest in the limit of a colored dead
leave model in the more traditional sense of finite-dimensional distributions. Since a dead leaves model $M$ is a random
tessellation (composed by its visible parts $(V_i)$, see Proposition~\ref{prop:deadleaves}), we may
use the results of Section \ref{colored}.  

\begin{definition}\label{def:colDeadLeave}
Let $\{C_0(\xx),\,\xx\in\R^2\}$ be a random field.
We denote by $I(C_0,r_0,r_1)$ the colored dead leaves model obtained from the
random tessellation $M(r_0,r_1)$ (see Definition \ref{def:CRT}).
If $C_0$ is the constant random filed with uniform marginals, that is, for all $\xun,\ldots,\xn$
and for all $c_1,\ldots,c_n\in(0,1)^n$ 
$$ 
P(C_0(\xun)\leq c_1,\ldots,C_0(\xn)\leq c_n)= 
P(C_0(\xun)\leq \min(c_1,\ldots,c_n))=\min(c_1,\ldots,c_n),
$$
we simply denote the colored dead leaves model by $I(r_0,r_1)$. 
In other words $I(r_0,r_1)$ is
obtained from the dead leaves model by independently coloring each leaf with a uniform distribution. 
\end{definition}

\begin{remark}
Note that the
definition of $I$ requires the knowledge of $(V_i)$, and not only the distribution of $\partial M$,  
since a visible part is not necessarily connected.
\end{remark}

\begin{remark}
Observe that, if $C_0$ is a stationary process, then the same is true for $I(C_0,r_0,r_1)$.
For instance, $I(r_0,r_1)$ is stationary.
\end{remark}

We now investigate the existence of a continuous prolongation of $I(C_0,r_0,r_1)$ at $r_0=0$. 
As in Section~\ref{sec:limit-boundary-set}, we assume that $1<\alpha\leq3$ and $r_1>0$ or that $\alpha>3$ and
$r_1\in(0,\infty]$. Recall that in that case the model $M(0,r_1)$ is not correctly defined although $M(r_0,r_1)$ is for
all $r_0\in(0,r_1)$. We have already taken interest
in the limit of this model as $r_0\to0$ in the sense of limit distribution of its boundary set (see
Section~\ref{sec:limit-boundary-set}) giving raise in degenerate limits even when the first order probability function
continuous prolongation $p(0,r_1,\cdot)$ is non-trivial. We now investigate the limit a colored dead leaves process defined
from $M(r_0,r_1)$ and a given color random field $C_0$. From Theorem~\ref{thm:limitProc}, for all random field $C_0$, it
is sufficient to study the limit of the corresponding second order partition process which we denote by
$\{R^{(2)}(r_0,r_1,\xx,\yy)\,:\,\xx,\yy\in\R^2\}$. The following remark and assumption will highly simplify this problem.

\begin{remark}\label{rem:zeroboundary} Applying Proposition~\ref{QKprop}, we have, for every  $\xx\in\R^2$,
  $Pr(\xx\in\partial M)=1-Q(\{\xx\})=1-E\nu(\inter(X))/E(\nu(X))$. This is clearly zero if we assume
\begin{hyp}{nullBoundary}
 $\nu(\partial Y)=0$ almost surely.
\end{hyp}
This is always true for \textit{non-degenerate} cases since $\partial Y$ typically is a curve. 
Note for instance that \refhyp{nullBoundary} holds under the assumption of Corollary~\ref{cor:assumptions}. 
In the sequel we assume \refhyp{nullBoundary} for convenience. 
\end{remark}

We obtain the following result.

\begin{proposition}\label{prop:fidiLimCDLM}
There exists a random process $I(C_0,0,r_1)$ such that 
$$
I(C_0,r_0,r_1)\fidi I(C_0,0,r_1)\mbox{ as }r_0\to0.
$$ 
\end{proposition}
\paragraph{proof}
Let $r_1>0$. In this proof, for all $r_0\in(0,r_1)$, we denote by
$\{V_i^{r_0}\}_i$ the visible parts of the dead leaves model $M(r_0,r_1)$. Using Theorem \ref{thm:limitProc}, it is enough to prove that there exists a
random process $R^{(2)}(0,r_1)$ such that $R^{(2)}(r_0,r_1)\fidi R^{(2)}(0,r_1)$ as
$r_0 \to 0$. More generally we may show that, for all $n\in\NN$, for all compact sets $K_1,\ldots,K_n$, 
the joint distribution of $\left( \indi(\exists i_j,\,K_{j}\subset V_{i_j}^{r_0}) \right)_{j=1}^n\in\{0,1\}^n$ converges as $r_0$
tend to the origin. Then applying
this to compact sets composed of two points gives the result. 
Let us denote in this proof section, for all $n\geq1$ and for all $\epsilon_1,\ldots,\epsilon_{n}\in\{0,1\}^{n}$, 
$$
p^{(n)}(\epsilon_1,\ldots,\epsilon_n):=
Pr(\left( \indi(\exists i_j,\,K_{j}\subset V_{i_j}^{r_0}) \right)_{j=1}^{n}=(\epsilon_j)_{j=1}^{n}.
$$
Observe now that, for all $n\geq2$ and for all $\epsilon_1,\ldots,\epsilon_{n-1}\in\{0,1\}^{n-1}$, 
$$
p^{(n)}(\epsilon_1,\ldots,\epsilon_{n-1},0)+p^{(n)}(\epsilon_1,\ldots,\epsilon_{n-1},1)
=p^{(n)}(\epsilon_1,\ldots,\epsilon_{n-1}).
$$
A simple induction on $n$ thus gives that it is sufficient to show
that the probability
$p^{(n)}(1,\ldots,1)=Pr(\exists i_1,\dots,i_n / K_1 \subset V_{i_1}^{r_0},\dots,K_n \subset V_{i_n}^{r_0})$
converges in $[0,1]$ as $r_0$ tends to 0. 
We may also assume without loss of generality that each $K_j$ contains at least two
distinct points. Otherwise $K_j$ is included in the interior of a visible part with probability one (see
Remark~\ref{rem:zeroboundary}). 
Finally, we claim that it is now enough to prove the convergence of  
$$
Pr(\exists i_1,\dots,i_n / K_1 \subset \inter{V_{i_1}^{r_0}},\dots,K_n \subset
\inter{V_{i_n}^{r_o}} \mbox{ and } t_{i_1}<t_{i_2}<\dots<t_{i_n}).
$$ 
This follows from the fact that the union of two compact sets is a compact set,
so that we may restrict ourselves to disjoint visible parts, by an elementary
induction. The probability above is defined as $Q^{(n)}_{r_0}(K_1,\ldots,K_n)$ as in~(\ref{eq:ncompacts}). 
After Proposition \ref{prop:ncompacts}, we know that
$$
 Q^{(n)}_{r_0}(K_1,\dots,K_n)=F_{r_0}^{(n)}(K_1,\dots,K_n)/G_{r_0}^{(n)}(K_1,\dots,K_n),
$$
with obvious notations.
Pick a $K_j$ and let $\delta$ denote its diameter. Recall that we have assumed $\delta>0$.
From~\refhyp{propy}, we have $\nu(RY\ominus \check  K_j)=0$ for all $R$ such that $2Ra_2$ is larger than $\delta$.
This implies that, 
$$
\eta(r_0,r_1)E\nu\left((\inter{X}\ominus \check  K_j)\cap  \Comp{(X\oplus\underline{\check K}_{j-1})}\right)
=\int_{r_0}^{r_1}
E_Y\nu\left((r\inter{Y}\ominus \check  K_j)\cap \Comp{(rY\oplus\underline{\check K}_{j-1})}\right)\,r^{-\alpha}\,dr 
$$
stays constant as soon as $r_0$ goes below $\delta/(2a_2)$. Here, for $j\geq1$, $\underline{\check K}_{j}$ is defined
in~(\ref{eq:defunderK}),  the case $j=1$ being obtained 
with the convention $\underline{\check K}_{0}=\emptyset$.
Hence, from~(\ref{eq:ncompacts1}), it is clear that, for
$r_0$ small enough, $F_{r_0}^{(n)}(K_1,\dots,K_n) \eta(r_0,r_1)^{-n}$ does not depend on $r_0$, with $\eta$ as in
(\ref{eq:defDensGeo}). 
On the other hand, for all $j=1,\ldots,n$, we have 
$$
E\nu(\check X \oplus \underline{\check K}_{j})
=\eta(r_0,r_1)\int_{r_0}^{r_1} r^{-\alpha}
E_Y \nu\left(rY \oplus\underline{\check K}_{j}\right) \,dr,
$$
where we recall that $E_Y$ is the expectation with respect to $Y$. Since the integrand is positive,
$\eta(r_0,r_1)^{-n}G_{r_0}^{(n)}(K_1,\dots,K_n)$ has a limit in $(0,\infty]$ (it is non zero since the $K_i$'s are non empty).
Simplifying by $\eta(r_0,r_1)^{-n}$ in the ratio defining $Q^{(n)}_{r_0}$ above, we obtain that it has a limit as
$r_0$ tend to the origin, which, as we claimed, is sufficient for showing Proposition~\ref{prop:fidiLimCDLM}. \cqfd

Note that in this proposition, we did not make any assumption on the value of
$\alpha$. However, we will see in Proposition \ref{prop:MeasModif} that this
result is interesting for $1<\alpha<3$, in which case there exist a measurable
version of our limit process. If $\alpha\geq3$, the limit process
simply is white noise.

We conclude this section by a simple corollary of Proposition~\ref{prop:fidiLimCDLM}, where we compute
the two-dimensional distributions of $I(C_0,0,r_1)$. We have, for all $\xx,\yy\in\R^2$, 
$$
Pr(R^{(2)}(r_0,r_1,\xx,\yy)=1)=p(r_0,r_1,\yy-\xx),
$$
where $p(r_0,r_1,\cdot)$ is defined in Section~\ref{sec:definitions} and computed in~(\ref{eq:Pgamma}). For all
$\xx\in\R^2$, $p(r_0,r_1,\xx)$ has a limit $p(0,r_1,\xx)$ and, from Lemma~\ref{lem:finitedistr} and
Remark~\ref{rem:bidim}, we easily obtain the two-dimensional distributions of $I(C_0,r_0,r_1)$. 
\begin{corollary}\label{cor:bidim}
For all $\xx,\yy\in\R^2$ and for all $r_0\geq0$, $(I(C_0,r_0,r_1,\xx),I(C_0,r_0,r_1,\yy))$ is a mixture of the two
(two-dimensional) random variables $(C_0(\xx),C_0(\yy))$ and $(C_0(\xx),C_1(\yy))$, where $C_1$ is an independent copy of
$C_0$, with respective weights $p(r_0,r_1,\yy-\xx)$ and $1-p(r_0,r_1,\yy-\xx)$.
\end{corollary}

Note furthermore that, as in the isotropic case~(\ref{eq:prolong_p_en0}), $\xx\mapsto p(0,r_1,\xx)$ is easily obtained by
taking the limit as $r_0\to0$
in the integrals of the RHS of~(\ref{eq:Pgamma}): the one in the numerator is always finite
by~(\ref{eq:gammatildeinfty}) and the one in the denominator is either positive or infinite by applying
Lemma~\ref{lem:gammatilde}. More precisely, we have, for all $\xx\in\R^2$, 
\begin{align}
\label{eq:p0degen}
&p(0,r_1,\xx)=\indi(\xx=0)\mbox{ if }\alpha\geq3,\\
\label{eq:p0}
&p(0,r_1,\xx)=
\frac{
\int_{0}^{r_1} \tilde{\gamma}(\xx/u) \, u^{2-\alpha}\,du}
{\int_{0}^{r_1} (2\tilde{\gamma}(0) - \tilde{\gamma}(\xx/u))\, u^{2-\alpha}\,du}
\mbox{ if }\alpha\in(1,3).
\end{align}
By Corollary~\ref{cor:bidim},~(\ref{eq:p0degen}) implies that $I(C_0,0,r_1)$ has independent 
(identically distributed if $C_0$ is stationary) samples if $\alpha\geq3$. 
For $1<\alpha<3$, (\ref{eq:p0}) and Proposition~\ref{prop:asympP} give that  $\xx\mapsto p(0,r_1,\xx)$ is a continuous
$\R^2\to[0,1]$ mapping.

\subsection{Preliminary properties of the limit process}

We have seen that, for $\alpha\geq3$, $I(0,r_1)$ is a white noise random field. On the contrary, for $\alpha\in(2,3)$, Proposition~\ref{prop:asympP}(\ref{item:asym2}) shows that the bi-dimensional
distributions given in Corollary~\ref{cor:bidim} (taking the constant field for $C_0$) exhibit interesting scaling
properties. We have so far only been interested in finite-dimensional distribution of the colored dead leaves model. 
Let us investigate how the simple scaling properties of the two-dimensional distributions influence the sample paths
properties of the model. The first property we may check is the existence of a measurable version of $I(C_0,0,r_1)$. 
A random field $\{Z(\omega,\xx)\,:\,\xx\in \R^2\}$ defined on $(\Omega,\G,P)$ is measurable if 
$(\omega,\xx)\mapsto Z(\omega,\xx)$ is a $(\Omega\times\R^2,\G\otimes\calB(\R^2))\to(\R,\calB(\R))$ (jointly)
measurable mapping (see e.g. \cite[Section 9.4]{samoro94}). Recall also that the random field $\{Z(\xx) ,\,\xx\in\R^2\}$ is
said to be stochastically continuous if, for all $\xx\in\R^2$, $Z(\yy)\cp Z(\xx)$ ($Z(\yy)$ converges to $Z(\xx)$ in
probability) as $\yy\to\xx$.

\begin{proposition}\label{prop:MeasModif}
Take $\alpha\in(1,3)$ and $r_1<\infty$.
Assume that $C_0$ is stochastically continuous. 
Then $I(C_0,0,r_1)$ is stochastically continuous.
If moreover $C_0$ has a measurable version or, equivalently, if, 
for all $G,H\in\calB(\R)$ and for all $\xx\in\R^2$, $\yy\mapsto Pr(C_0(\xx)\in G,\,C_0(\yy)\in H)$ is
  a $(\R^2,\calB(\R^2))\to([0,1],\calB([0,1]))$ measurable mapping, then 
there exists a measurable modification of
$I(C_0,0,r_1)$, that is a measurable process
$\tilde{I}(C_0,0,r_1)$ defined on the same probability space such that, for all $\xx\in\R^2$,
$Pr(\tilde{I}(C_0,0,r_1,\xx)=I(C_0,0,r_1,\xx))=1$.
\end{proposition}
\paragraph{proof}
For convenience we write $I$ for $I(C_0,0,r_1)$ in this proof. 
The two-dimensional distributions of $I$ are given in Corollary~\ref{cor:bidim}. We have seen that $\xx\mapsto
p(0,r_1,\xx)$ defined by~(\ref{eq:p0}) is a continuous mapping. We obtain, for all $\xx,\yy\in\R^2$
and for all $\epsilon>0$, 
$$
Pr(|I(\xx)-I(\yy)|>\epsilon)
\leq Pr(|C_0(\xx)-C_0(\yy)|>\epsilon)+(1-p(0,r_1,\yy-\xx)),
$$
which tends to zero as $\yy\to\xx$ for $C_0$ stochastically continuous and since $p(0,r_1,0)=1$.
Hence the first part of the proposition. For the second part, we apply \cite[Theorem 9.4.2]{samoro94}. We have to check
two conditions on $I$, namely 
\begin{enumerate}[{\rm (i)}]
\item \label{item:C1} There exists a countable set $S\subset\R^2$ such that for all $\xx\in\R^2$, there exists a
  sequence $(\xx_k)$ such that $I(\xx_k)\cp I(\xx)$ as $k\to\infty$ and $\xx_k\in S$ for all $k$. 
\item \label{item:C2} For all $G,H\in\calB(\R)$ and for all $\xx\in\R^2$, $\yy\mapsto Pr(I(\xx)\in G,\,I(\yy)\in H)$ is
  a $(\R^2,\calB(\R^2))\to([0,1],\calB([0,1]))$ measurable mapping.
\end{enumerate}
Condition~(\ref{item:C1}) is a consequence of stochastic continuity. 
Moreover, condition~(\ref{item:C2}) must be satisfied by the random field $C_0$ if it has a measurable version (and
since $C_0$ satisfies~(\ref{item:C1}) as a consequence of stochastic continuity, condition~(\ref{item:C2}) on $C_0$
implies the existence of a measurable version for $C_0$). Condition~(\ref{item:C2}) is easily checked on $I$ by using
the two-dimensional distributions given in Corollary~\ref{cor:bidim}.\cqfd

\begin{corollary}\label{cor:constantunifcolor}
The colored dead leave model $I(0,r_1)$ is stochastically continuous and it admits a measurable modification if and only
if $\alpha\in(1,3)$.
\end{corollary}
\paragraph{proof}
For $\alpha\geq3$,  $I(0,r_1)$ is a white noise and thus does not have a measurable version
(see \cite[Example 9.4.3]{samoro94}) and neither is stochastically continuous.
For $\alpha\in(1,3)$ the constant random field trivially satisfies the assumptions of Proposition~\ref{prop:MeasModif}.
\cqfd

Let us conclude this section by noting that this result is coherent with a result of D. Mumford and B. Gidas which, in a paper investigating stochastic
models for natural images, \cite{MumfordGidas}, show that a non-constant scale-invariant twodimensional random
field cannot be a random function, and has to be a random distribution
(generalized function). In the framework of the dead leaves model, this
corresponds to the case $\alpha=3$,
where there is no measurable version of the model. Indeed, the
scale invariant case would correspond to $\alpha=3$ when $r_0 \to 0$ and $r_1
\to \infty$ simultaneously. In this case the limit of the two-dimensional distribution $(I(\xx,r_0,r_1),I(\yy,r_0,r_1))$
depends on the ratio  between $\log r_0$ and $\log r_1$, but never on $(\xx,\yy)$ (see Proposition
\ref{deg_limits_synthese} (\ref{eq:deg3})) and we find a
measurable limit only in the case where this limit is a constant field. In the next section, we
investigate the regularity of our limit model in the framework of Besov spaces.

\begin{figure}
\begin{center}
\includegraphics[width=7cm]{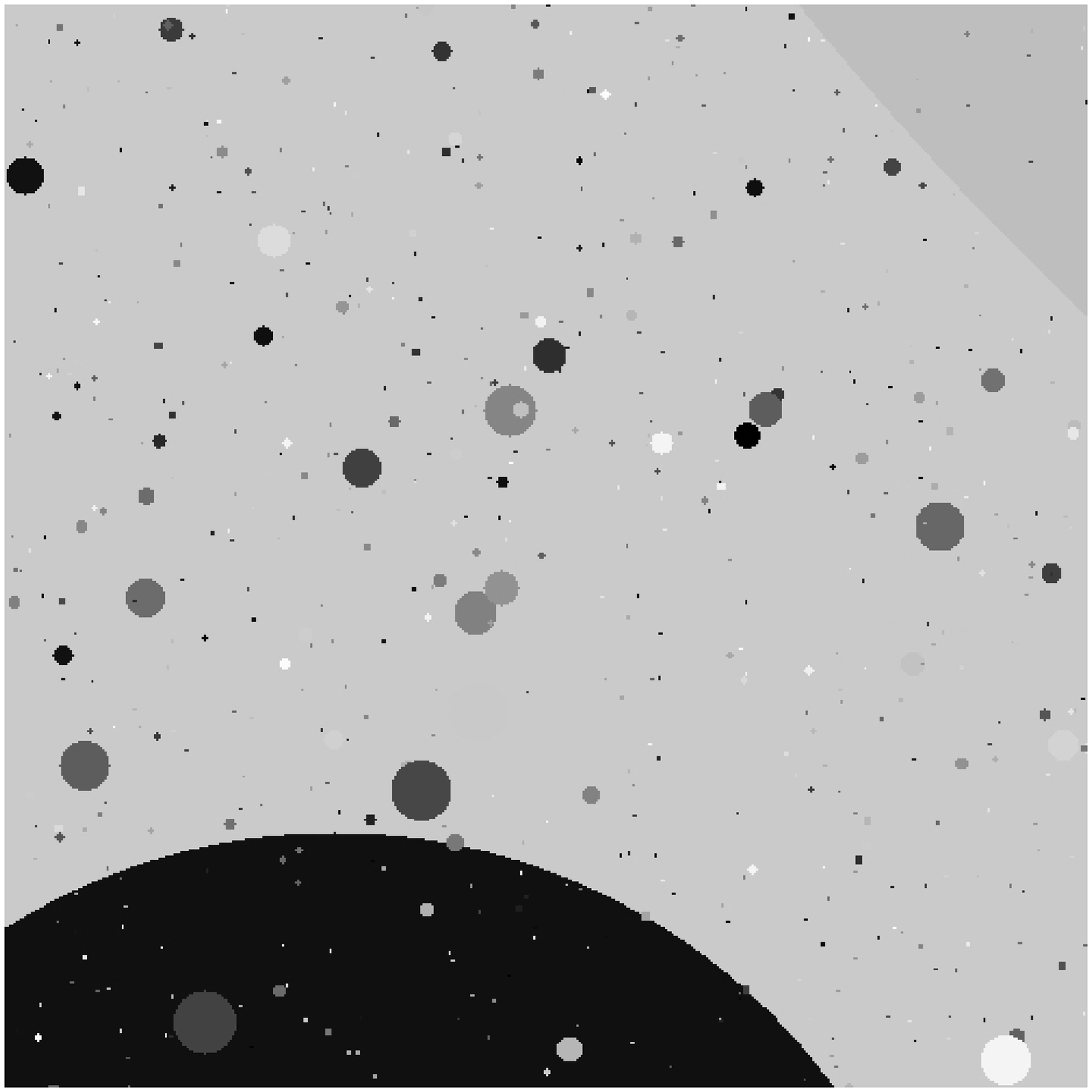} \qquad
\includegraphics[width=7cm]{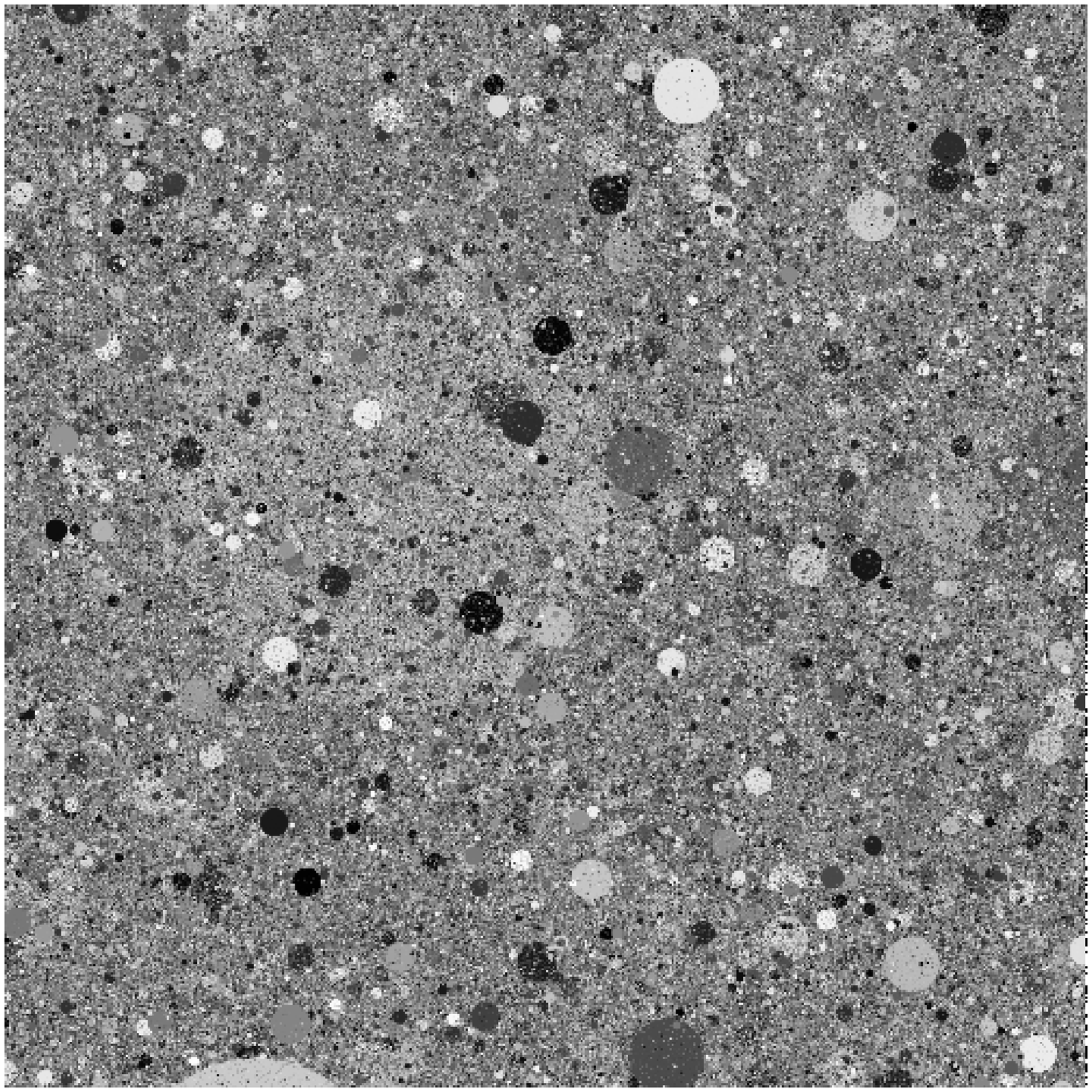}

\caption{Illustration of the degenerate cases of Proposition \ref{deg_limits_synthese}. Left: case
  (\ref{eq:deg1}), the process converges to a constant (random) function. Right: case
  (\ref{eq:deg2}), the process converges to a white noise.}
\label{fig:degenerate}
\end{center}
\end{figure}

\begin{figure}
\begin{center}
\includegraphics[width=7cm]{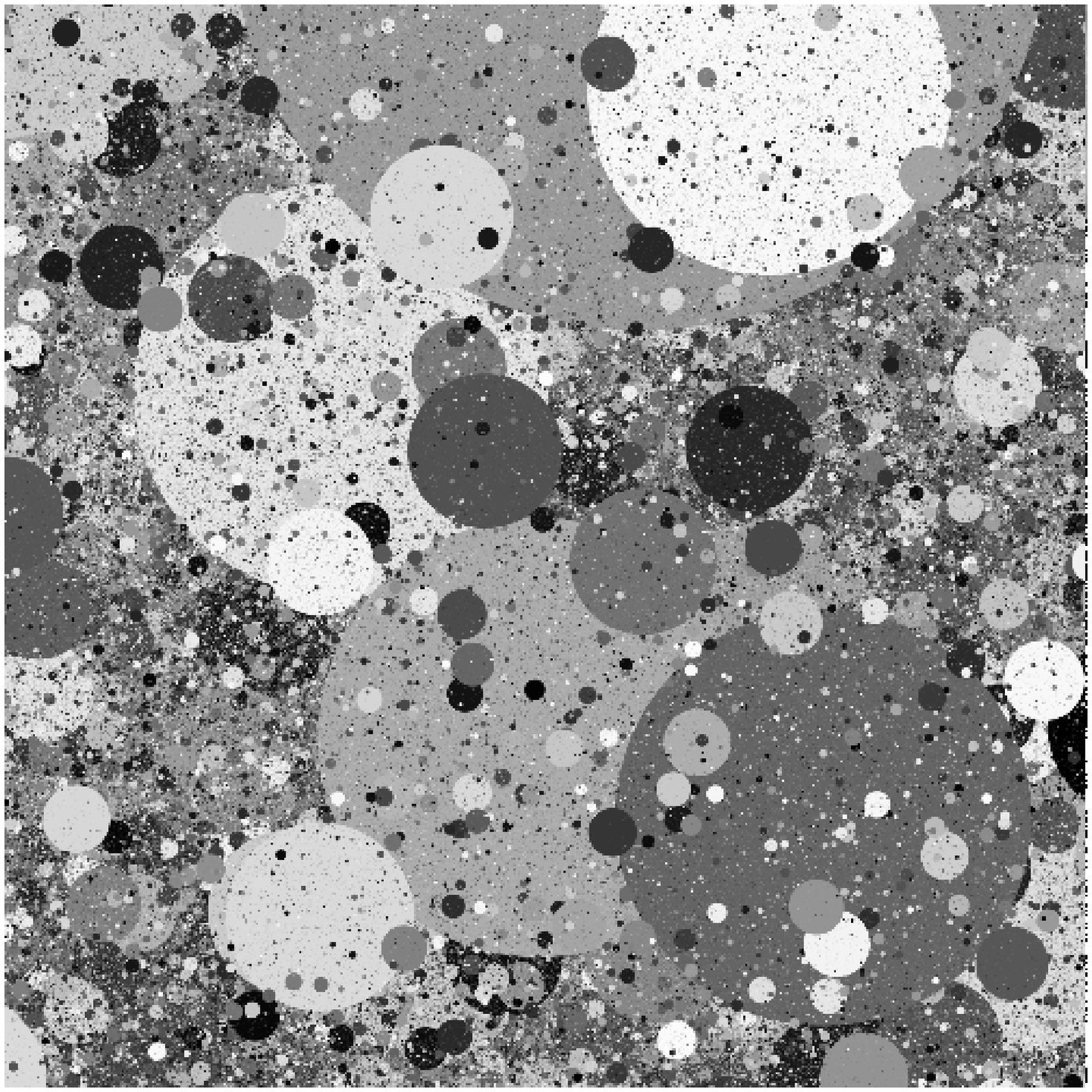} \qquad
\includegraphics[width=7cm]{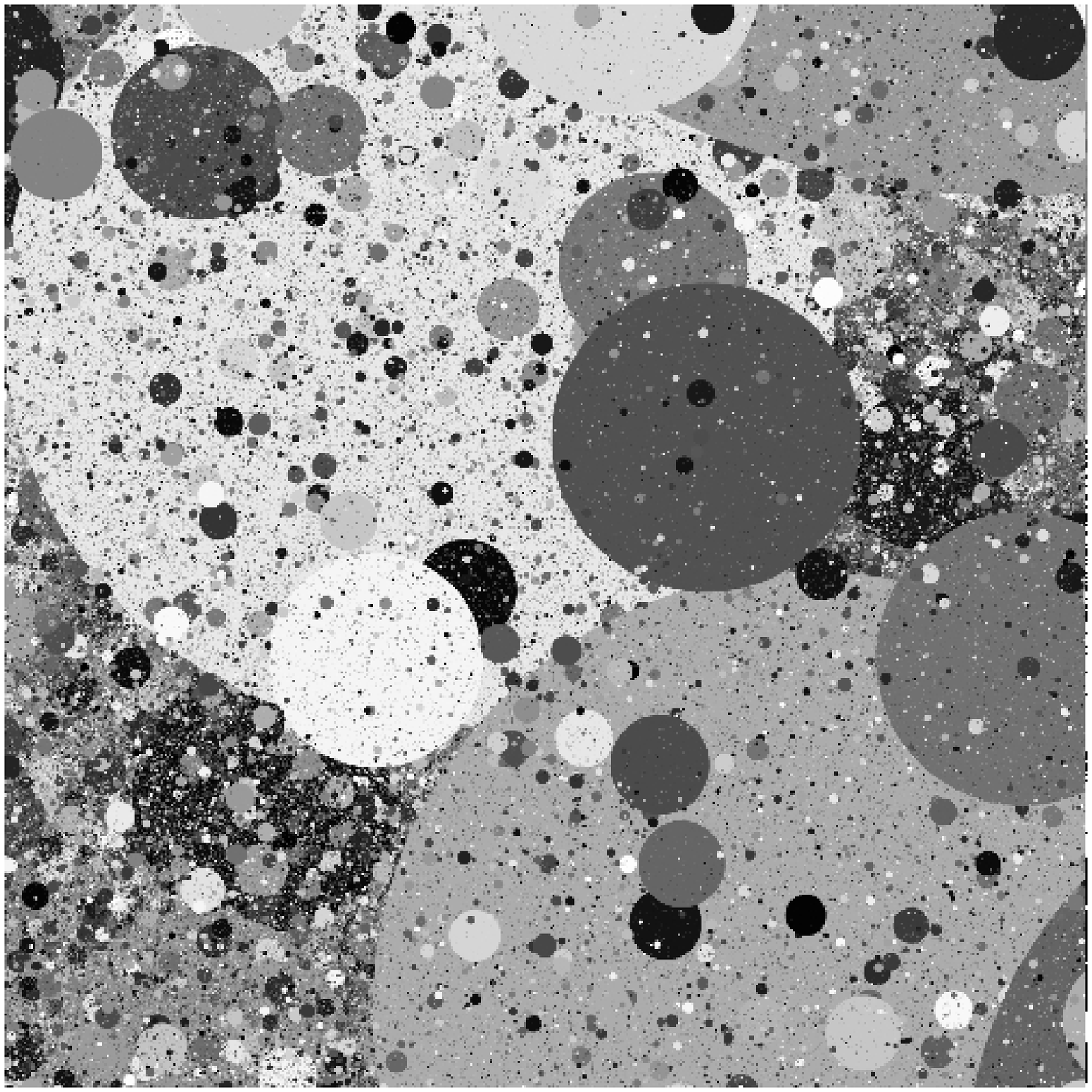}
\vskip .8cm
\includegraphics[width=7cm]{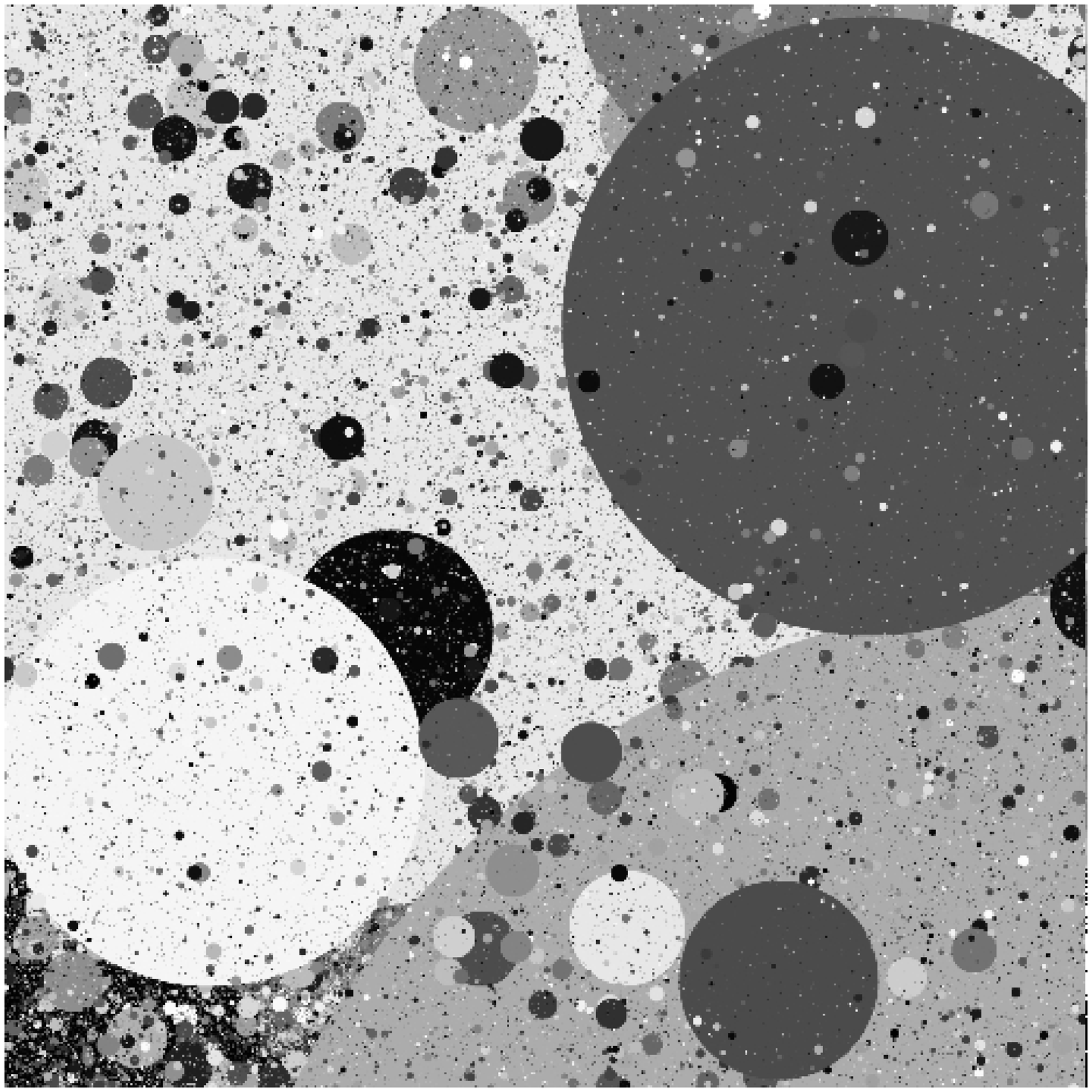}  \qquad
\includegraphics[width=7cm]{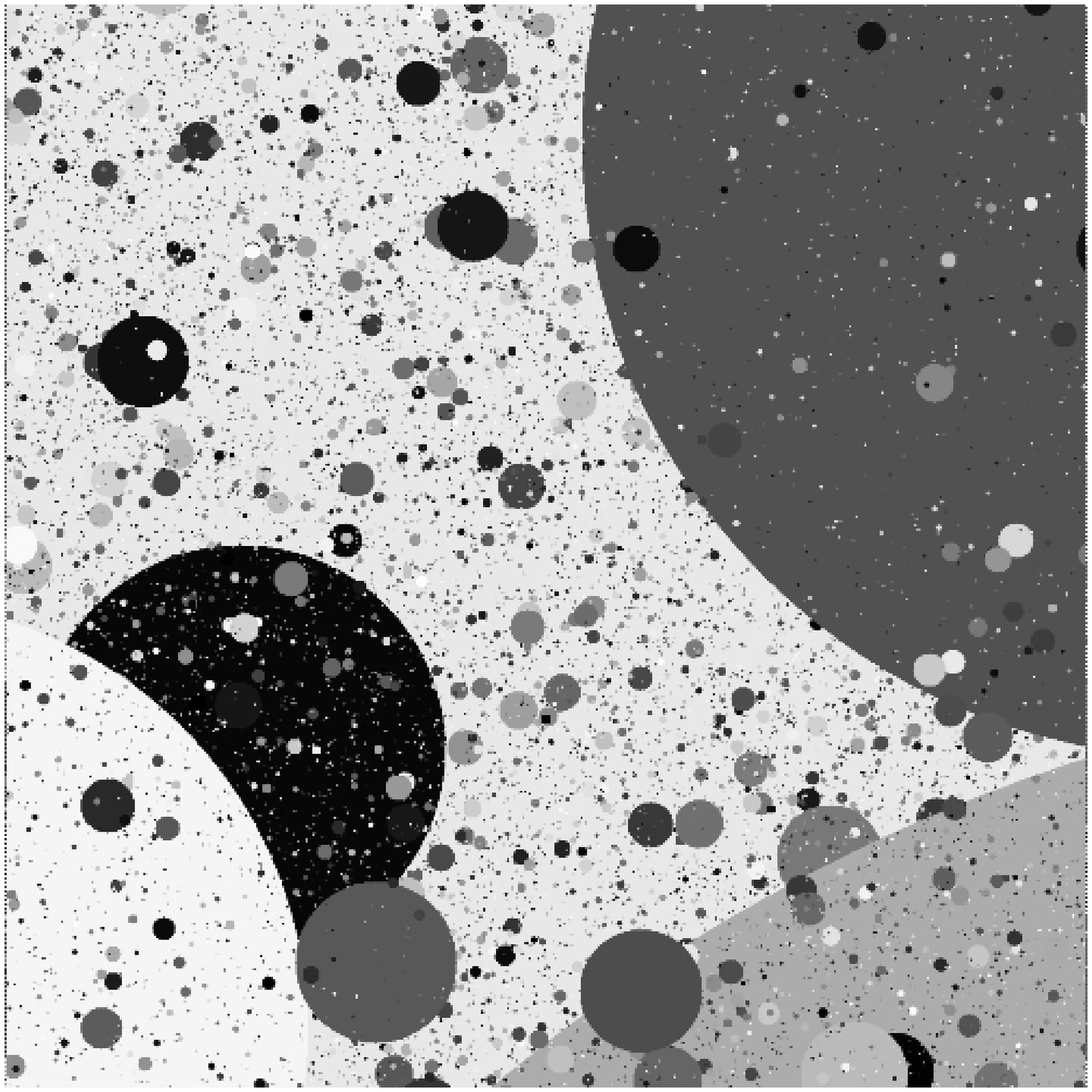}

\caption{Illustration of the convergence of the model when
  $\alpha=2.9$ and $r_0 \to 0$. The first image (up left) is a realization when $r_0=1$,
  $r_1=10000$, on a window of size $10000 \times 10000$. The three next images
  are zooms on this realization, the zoom factor being two each time.}
\label{fig:zooms}
\end{center}
\end{figure}

Before proceeding, we illustrate the properties of the limit model with some
simulations. In figure \ref{fig:degenerate}, we show two examples illustrating
Proposition \ref{deg_limits_synthese}. Images are simulated using the perfect
simulation methods explained in Section \ref{sec:deadleaves}, and gray levels
are uniform between 0 and 255. In the first example (left) we illustrate point
(\ref{eq:deg1}); 
$\alpha=2.5$, and $r_1 \to \infty$. The image is
of size $1000 \times 1000$, $r_0=1$, $r_1=100000$; the process converges to a
constant function. In the second example (right), we illustrate point (\ref{eq:deg2});
$\alpha=3.5$, and $r_0 \to 0$. The image is of
size $10000 \times 10000$, $r_0=1$, $r_1=10000$, the process converges to
white noise. In figure \ref{fig:zooms} we illustrate the convergence of
$I(r_0,r_1)$ when $r_0 \to 0$ and $\alpha=2.9$. The first image is of size
$10000 \times 10000$, $r_0=1$, $r_1=10000$. The next three images are zooms of
the same realization of the model (the zoom factor is two each time).

\subsection{Some sample paths properties} 

If $r_0>0$, the process $I(C_0,r_0,r_1)$ have paths for which occlusions influence the smoothness in a rather simple way: it
introduces discontinuities along $\partial M(r_0,r_1)$ and, between these discontinuities, the smoothness is driven by the
properties of $C_0$. For instance, in the case of a constant $C_0$, $I(r_0,r_1)$ simply is a piecewise constant
process, the pieces being connected components of $\R^2\backslash\partial M(r_0,r_1)$. In the simple case where $\partial Y$
has finite length, $\partial M(r_0,r_1)$ has locally finite length almost
surely (as an easy consequence of Lemma \ref{lemma:reduit}) and  thus, under suitable regularity assumption for $C_0$, is locally
of bounded variation. This is contradicted by empirical experiments. In \cite{Gousseau01}, by
investigating the distribution of sizes of "bilevel" sets in natural images
(up to the smallest available scale),
it is shown that this smoothness model is erroneous. In
practice this means that a denoising approach relying on the \textit{a priori}
of a piecewise smooth image may interpret
small objects as noise and, therefore, may result as a non-negligible loss of
information. This is well known in the field of image processing. For
instance, variational methods in the space of functions with bounded
variation, such as the famous
Rudin-Osher denoising scheme, \cite{RudinOsherFatemi}, are known to erase textured area.
The efficiency of these methods is dued to their ability to preserve the large
scale geometric structure of images, but they often fail in handling the small
scales structure of the image. A recent approach to overcome this difficulty
has been proposed by Y. Meyer, see \cite{MeyerBV}, introducing a new
functional space to account for textured regions in images, see also
\cite{VeseBV}, \cite{AubertBV}.

The idea here is to propose an \textit{a priori} model which includes the existence of small objects. How this model
could be used efficiently in practice, say for denoising, may be a difficult question. In non-parametric regression (including
the denoising problem), an important issue is to model the smoothness of the target function in terms
of standard smoothness classes and, in many cases, this issue is closely related to an \textit{a priori} statistical model
(see e.g. \cite{abram98}). The object of this section is to provide some insights in the smoothness properties of the model
$I(0,r_1)$ which we believe to be a promising \textit{a priori} statistical model for natural images. We will consider Besov
smoothness spaces, which have been standard spaces in nonparametric regression since it was promoted in the celebrated work
of Donoho and Johnstone (see e.g.  \cite{donohojohnstone98}).
Of course not every value of $\alpha$ is relevant for the model $I(0,r_1)$. We will avoid $\alpha\geq3$ be cause in this case
$I(0,r_1)$ simply is a white noise (see
Corollary~\ref{cor:constantunifcolor}). In fact, the most suitable range for the
modeling natural images by $I(0,r_1)$ is $2<\alpha<3$. Indeed, in this case $p(0,r_1,\cdot)$ exhibit interesting scaling
properties at small scales (see Proposition~\ref{prop:asympP}) which correspond to
the empirical observations (see \cite{Ruderman96}). Hence in the
following we will only consider this case when $r_0=0$. When $r_0>0$, we will consider the usual case $\alpha\in(1,3]$ and
$0<r_0<r_1<\infty$ or $\alpha>3$ and $0<r_0<r_1\leq\infty$.

Since in all these cases, there exists a $(\R^2,\calB(\R^2))\to(\R,\calB(\R))$ measurable 
modification $\xx\mapsto\tilde{I}(\omega,\xx)$ of $I(r_0,r_1)$, it makes sense to investigate whether $\tilde{I}$ belongs to some
given functional spaces. We only take interest in local smoothness. Since the process is stationary, we may consider
its restriction to the cube $[0,1]^2$ without loss of generality. We will take
interest in Besov spaces and, as a byproduct, obtain that our limit process $I(0,r_1)$ is not of bounded variation, which is coherent with the results
of \cite{Gousseau01} in view of natural images modeling. Let
$s\in(0,1)$, $p\in[1,\infty]$ and $q\in[1,\infty]$. The
Besov space $B_p^{s,q}([0,1]^2)$  (see e.g. \cite{devore91} or \cite{meyer90}) is the Banach space endowed with the following
norm 
$$
|f|_{B_p^{s,q}}:=|f|_p+\left(\int_{u>0}\left(\omega(f,u)_p \, u^{-{s}}\right)^q\,\frac{du}{u}\right)^{\frac1q},
$$
where $|\cdot|_p$ is the usual $L^p([0,1]^2)$ norm and $\omega(f,u)_p$ is the  $L^p([0,1]^2)$ modulus of smoothness of
$f$ at scale $u$, that is
$\omega(f,u)_p:=\sup_{|\yy|<u}|\Delta(f,\yy)|_p$,
where $\Delta(f,\yy)$ is the difference operator applied to $f$ with step $\yy$ on $[0,1]^2$, that is, the
mapping  $\xx\mapsto (f(\xx+\yy)-f(\xx))\indi(\xx\in[0,1]^2,\,\xx+\yy\in[0,1]^2)$.

Finally let us note that the results below can be easily generalized to 
a colored dead leaves process $I(C_0,r_0,r_1)$ with $C_0$ satisfying the assumptions of
Proposition~\ref{prop:MeasModif}, but, in general, they would depend on $C_0$. Here we focus on the basic
properties implied by the geometrical construction of the model. Therefore we only consider the case of $I(r_0,r_1)$,
that is with $C_0$ being the constant field, because then the variations of the colors are directly connected to the geometry
of the model. For similar reasons, if $\partial Y$ has some irregularity, it may influence the smoothness of
$I(r_0,r_1)$. This is clear from the example given in the end of Section~\ref{sec:definitions}. Hence the additional
assumptions \refhyp{gamma1} or \refhyp{gamma2} in the following results.  

\begin{proposition}\label{prop:BesovEsp}
For all $p\in[1,\infty]$ and for all $s\in(0,1)$, the following assertions hold true. 
\begin{enumerate}[{\rm (a)}]
\item \label{item:Cas1} Let $1<\alpha\leq3$ and $0<r_0<r_1<\infty$ or if $\alpha>3$ and
  $0<r_0<r_1\leq\infty$. Under~\refhyp{gamma2}, 
\begin{equation}\label{eq:besovCas1} 
E\left[|\tilde{I}|_{B_p^{s,p}}^p\right] < \infty \Leftrightarrow s<1/p.
\end{equation}
\item \label{item:Cas2} Let $2<\alpha<3$ and $0=r_0<r_1<\infty$. Under~\refhyp{gamma1}, 
\begin{equation}\label{eq:besovCas2} 
 E\left[|\tilde{I}|_{B_p^{s,p}}^p\right] < \infty \Leftrightarrow s<\frac{3-\alpha}p.
\end{equation}
\end{enumerate}
\end{proposition}

\paragraph{proof}
In this proof we write $A\asymp B$ if there exists a constant $c$ (possibly depending on the constants $s,p,\alpha,r_0,r_1$
and $\tilde{\gamma}$) such that $\frac1c B \leq A \leq c B$. It is more convenient to use the modified modulus of
smoothness   
$$
w(f,u)_p:=(1\vee u^{-2})\int_{|\yy|<u}|\Delta(f,\yy)|_p\,d\yy.
$$
which satisfies $w(\cdot,\cdot)_p\asymp \omega(\cdot,\cdot)_p$  (see~\cite{devore91}).
Because $\tilde{I}$ is measurable, we may use the Fubini Theorem. We obtain
\begin{align}
\nonumber 
E[|\tilde{I}|_{B_p^{s,p}}^p]&\asymp  
E[|\tilde{I}|_p]^p + \int_{u>0}\int_{|\yy|<u}E[|\Delta(\tilde{I},\yy)|_p^p]\,(1\vee 
u^{-2})u^{-p{s}}\,d\yy\,\frac{du}{u}\\  
\nonumber 
&\asymp 1+\int_{u>0} \int_{|\yy|<u\wedge1} (1-p(r_0,r_1,\yy)) \,(1\vee u^{-2})u^{-p{s}}\,d\yy\,\frac{du}{u}\\
\label{eq:EtildI}
& = 1+ 2 
\int_{u>0} \int_{|\yy|<u\wedge1}
\frac{\int_{r_0}^{r_1}(\tilde{\gamma}(0)-\tilde{\gamma}(\yy/v))\,v^{2-\alpha}\,dv}
{\int_{r_0}^{r_1}(2\tilde{\gamma}(0)-\tilde{\gamma}(\yy/v))\,v^{2-\alpha}\,dv} 
\,(1\vee u^{-2})u^{-p{s}}\,d\yy\,\frac{du}{u}. 
\end{align}
Assume we are in the case~(\ref{item:Cas1}). By Lemma~\ref{lem:gammatilde}, the denominator and the numerator of the fraction
in the RHS of~(\ref{eq:EtildI}) are continuous mapping of $\yy\in\R^2$. Moreover, since
$0\leq\tilde{\gamma}\leq\tilde{\gamma}(0)$ the denominator is bounded away from zero and infinity independently of $\yy$.
Under~\refhyp{gamma2}, the numerator $\asymp(1\wedge|\yy|)$.  Case~(\ref{item:Cas1}) is then easily achieved. 

We now consider case~(\ref{item:Cas2}).
By Lemma~\ref{lem:gammatilde} and since $\alpha<3$, the denominator writes
\begin{equation}\label{eq:EtildINumerator}
\int_0^{r_1}(2\tilde{\gamma}(0)-\tilde{\gamma}(\yy/v))\,v^{2-\alpha}\,dv\asymp 
\tilde{\gamma}(0)r_1^{3-\alpha}/(3-\alpha)>0.  
\end{equation}
Concerning the numerator, a change of variable gives, for all $\yy\neq0$,
$$
\int_0^{r_1}(\tilde{\gamma}(0)-\tilde{\gamma}(\yy/v))\,v^{2-\alpha}\,dv
=|\yy|^{3-\alpha}\int_0^{r_1/|\yy|}(\tilde{\gamma}(0)-\tilde{\gamma}(\yy/(|\yy|v)))\,v^{2-\alpha}\,dv
$$
Consequently, since $\alpha<3$, under~\refhyp{gamma1},
$$
\sup_{|\bz|=1}\int_0^{\infty}(\tilde{\gamma}(0)-\tilde{\gamma}(\bz/v))\,v^{2-\alpha}\,dv<\infty.
$$ 
By~(\ref{eq:gammatildeinfty}), we have, for all $\bz$ such that $|\bz|=1$ and for all
$r\in(0,1/(2a_2))$, 
$$
\int_0^{r}(\tilde{\gamma}(0)-\tilde{\gamma}(\yy/(|\yy|v)))\,v^{2-\alpha}\,dv
= \int_0^{r} \tilde{\gamma}(0)\,v^{2-\alpha}\,dv= \frac{\tilde{\gamma}(0)}{3-\alpha}
r^{3-\alpha}. 
$$
By Lemma~\ref{lem:gammatilde}, we know that $\tilde{\gamma}(0)-\tilde{\gamma}(\yy/(|\yy|v))$ is non-negative, so that
the last three equations finally give
$$
\int_0^{r_1}(\tilde{\gamma}(0)-\tilde{\gamma}(\yy/v))\,v^{2-\alpha}\,dv
\asymp |\yy|^{3-\alpha} (|\yy|\vee1)^{\alpha-3}.
$$
From~(\ref{eq:EtildI}),~(\ref{eq:EtildINumerator}) and the last equations, we obtain
$$
E[|\tilde{I}|_{B_p^{s,p}}^p] \asymp
1+ 
\int_{u>0} \left(\int_{|\yy|<u\wedge1} |\yy|^{3-\alpha}\,d\yy\right)
\,(1\vee u^{-2})u^{-p{s}}\,\frac{du}{u}.
$$
Hence the result.
\cqfd

An almost sure smoothness result follows immediately.

\begin{corollary}\label{cor:BesovAS}
For all $p\in[1,\infty]$, for all $q\in[1,\infty]$ and for all $s\in(0,1)$, the following assertions hold true. 
\begin{enumerate}[{\rm (a)}]
\item \label{item:ASCas1} If $1<\alpha\leq3$ and $0<r_0<r_1<\infty$ or if $\alpha>3$ and $0<r_0<r_1\leq\infty$, under~\refhyp{gamma2},
\begin{equation}\label{eq:besovASCas1} 
s<1/p \Rightarrow \tilde{I}\in B_p^{s,q} \mbox{ a.s.}
\end{equation}
\item \label{item:ASCas2} If $2<\alpha<3$ and $0=r_0<r_1<\infty$, under~\refhyp{gamma1},
\begin{equation}\label{eq:besovASCas2} 
s<\frac{3-\alpha}p \Rightarrow \tilde{I}\in B_p^{s,q} \mbox{ a.s.}
\end{equation}
\end{enumerate}
\end{corollary}

\paragraph{proof}
We simply use that if a non-negative random variable has finite expectation, it is necessary finite almost surely. Then well
known inclusions of Besov spaces give the claimed results.
\cqfd

It is worthwhile to observe that if $f\in B_p^{s,p}$ in the parameter region $s>1/p$, then $f$ is continuous. Hence the
condition for Besov smoothness given in case~(\ref{item:Cas1}) is the almost (because the boundary $s=1/p$ is not included)
weakest condition for an almost surely discontinuous $\tilde{I}$. In other words the Besov smoothness
condition~(\ref{eq:besovCas1}) do not distinguish $\tilde{I}$ with any other standard discontinuous random field from a Besov
smoothness point of view. In particular it does not depend on $\alpha$. In contrast the case~(\ref{item:Cas2}) exhibits some
interesting behavior. This behavior is of course closely related to the asymptotical expansion found in Proposition~\ref{prop:asympP}
in the isotropic case. 
Moreover in this interesting case, as we anticipated before, the range of Besov spaces in which lies
  $I(0,r_1)$ is strictly bigger than the space $BV([0,1]^2)$. Indeed, it is well known that $BV([0,1]^2)
  \subset B_1^{\infty,1} \subset B_1^{s,1}$ for any $s>1$, and moreover that for any $f \in
  BV([0,1]^2)$, $||f||_{BV}\geq C||f||_{B_1^{\infty,1}}\geq C' ||f||_{B_1^{s,1}}.$
Therefore, as a consequence of Proposition \ref{prop:BesovEsp} (\ref{item:Cas2}), we have that,
for any $r_1>0$, $\alpha>2$, $E||I(0,r_1)||_{BV}=\infty$.

\paragraph{Concluding remarks}
In this paper we introduced two dimensional random fields obtained
as limits of sequences of dead leaves models. As advocated earlier, we
believe them to be promizing a priori models for natural images. First it
accounts for the geometry of natural images (in short the occlusion
phenomenon). Second, we used a distribution of objects sizes based on
empirical statistical properties of images, as promoted in
previous works .
Third, getting rid of small sizes ``cut off'' allows for
regularity results that are also compatible with empirical observations.

Now, how this model can be used in practice remains
open. On the one hand, the almost sure smoothness properties may be taken
as a functional a priori for images. For instance, in a denoising framework,
Corollary \ref{cor:BesovAS} provides a class of smoothness hypotheses only
driven by parameter $\alpha$. For more elaborate tasks such as shape
extraction, it is clear that parameter $\alpha$ is not sufficient, and that
geometrical properties of the model, depending on the distribution of objects
shapes, should be investigated. To this end further work has to be done on the
dead leaves model. On the other hand, the model could be used as a Bayes a
priori, but this would need tractable knowledge of the
dependance structure of the model.

\paragraph{Acknowledgements} We would like to thank Charles Bordenave for fruitful discussions about
Palm calculus. This work was partially supported by the CNRS working group GDR ISIS.

\bibliographystyle{alpha}
\bibliography{def,francois}

\end{document}